\documentclass[11.5pt,leqno]{article}
\usepackage{amsmath, amsthm, amsfonts, color, amssymb}
\usepackage{mathrsfs}
\usepackage{url}
\usepackage{hyperref}

\newcommand{\scr}[1]{\mathscr #1}

\newcommand{\black}{\color{black}}
\numberwithin{equation}{section} \theoremstyle{remark}

\title{{\bf  Exact Wasserstein Convergence Rate   for Conditional Empirical Measures of  Killed Diffusion Processes }\footnote{Supported in part by the National Natural Science Foundation of China (Grant Nos.~12531007 and 12401174), the National Key R\&D Program of China (Grant No.~2023YFA1010400), and the Alexander von Humboldt Foundation.} }
\author{{\bf    Feng-Yu Wang$^{a)}$, Bingyao Wu$^{b)}$}\\
\footnotesize{ $^{a)}$Center for Applied Mathematics and KL-AAGDM, Tianjin University, Tianjin 300072, China}\\
 \footnotesize{$^{b)}$School of Mathematics and Statistics,
 Fujian Normal University, Fuzhou 350007, China}\\
\footnotesize{wangfy@tju.edu.cn, bingyaowu@163.com.}}

\begin{document}
\allowdisplaybreaks
\maketitle
\def\ff{\frac} \def\ss{\sqrt}
\def\kk{\kappa} \def\vv{\varepsilon}
\def\dd{\delta} \def\DD{\Delta} \def\rr{\rho}
\def\gg{\gamma} \def\nn{\nabla} \def\pp{\partial}
\def\E{\mathbb E} \def\d{\text{\rm{d}}} \def\bb{\beta}
\def\si{\sigma} \def\beg{\begin} \def\beq{\begin{equation}}
\def\e{\text{\rm{e}}} \def\tt{\tilde} \def\P{\mathbb P}
\def\ll{\lambda} \def\W{\mathbb W} \def\C{\scr C}

\def\R{\mathbb R}  \def\ff{\frac} \def\ss{\sqrt} \def\BB{\mathbf
B}
\def\N{\mathbb N} \def\kk{\kappa} \def\m{{\bf m}}
\def\ee{\varepsilon}\def\ddd{D^*}
\def\dd{\delta} \def\DD{\Delta} \def\vv{\varepsilon} \def\rr{\rho}
\def\<{\langle} \def\>{\rangle}
  \def\nn{\nabla} \def\pp{\partial} \def\E{\mathbb E}
\def\d{\text{\rm{d}}} \def\bb{\beta} \def\aa{\alpha} \def\D{\scr D}
  \def\si{\sigma} \def\ess{\text{\rm{ess}}}\def\s{{\bf s}}
\def\beg{\begin} \def\beq{\begin{equation}}  \def\F{\scr F}
\def\Ric{\mathcal Ric} \def\BBess{\text{\rm{Hess}}}
\def\e{\text{\rm{e}}} \def\ua{\underline a} \def\OO{\Omega}  \def\oo{\omega}
 \def\tt{\tilde}\def\[{\lfloor} \def\]{\rfloor}
\def\cut{\text{\rm{cut}}} \def\P{\mathbb P} \def\ifn{I_n(f^{\bigotimes n})}
\def\C{\scr C}      \def\aaa{\mathbf{r}}     \def\r{r}
\def\gap{\text{\rm{gap}}} \def\prr{\pi_{{\bf m},\varrho}}  \def\r{\mathbf r}
\def\Z{\mathbb Z} \def\vrr{\varrho} \def\ll{\lambda}
\def\L{\scr L}\def\Tt{\tt} \def\TT{\tt}\def\II{\mathbb I}
\def\i{{\rm in}}\def\Sect{{\rm Sect}}
\def\M{\mathbb M}\def\Q{\mathbb Q} \def\texto{\text{o}} \def\LL{\Lambda}
\def\Rank{{\rm Rank}} \def\BB{\mathbb B} \def\i{{\rm i}} \def\BBR{\hat{\R}^d}
\def\to{\rightarrow}\def\l{\ell}\def\iint{\int}\def\gg{\gamma}
\def\EE{\scr E} \def\W{\mathbb W}
\def\A{\scr A} \def\Lip{{\rm Lip}}\def\S{\mathbb S}
\def\B{\scr B}\def\Ent{{\rm Ent}} \def\i{{\rm i}}\def\itparallel{{\it\parallel}}
\def\g{{\mathbf g}}\def\Sect{{\mathcal Sec}}\def\T{\mathcal T}\def\BBB{{\bf B}}
\def\f{\mathbf f} \def\g{\mathbf g}\def\BBL{{\bf L}}  \def\BBG{{\mathbb G}}
\def\BBd{{D^E}} \def\BBdP{D^E_\phi} \def\BBdd{{\bf \dd}} \def\BBs{{\bf s}} \def\GA{\scr A}
\def\BBg{{\bf g}}  \def\BBdd{\psi_B} \def\supp{{\rm supp}}\def\div{{\rm div}}
\def\ddiv{{\rm div}}\def\osc{{\bf osc}}\def\1{{\bf 1}}\def\BBD{\mathbb D}\def\GG{\Gamma}
\def\fdd{\stackrel{f.d.d.}{\longrightarrow}}\def\H{\mathbb H}

\begin{abstract}
Let $M$ be a $d$-dimensional compact connected Riemannian manifold $M$ with boundary $\pp M$, and let $X_t$ be the diffusion process generated by $L:=\DD+\nn V$ for some $V\in C^2(M)$ which is killed on the boundary. For any initial distribution $\nu$ which is not fully supported on $\pp M$,
let $\mu_t^\nu$ be the conditional  empirical measure for the  diffusion process  under the condition that the process is not yet killed at time $t$, and let $\mu_0$ be the associated
quasi-stationary distribution.  We derive the sharp convergence rate of $\mu_t^\nu$ to $\mu_0$ in the quadratic Wasserstein distance $\W_2$ as $t\to\infty$. In particular,    if  $\nu(B(x,r))\asymp r^k$ for some $k\in [0,d]$ and all $r\in [0,1]$ uniformly in $x\in {\rm supp}\nu$,
then  for large $t>0$
$$
\W_2(\mu_t^\nu,\mu_0)\asymp_\nu
\begin{cases}
t^{-1},&d-k<6,\\
t^{-1}\sqrt{\log t},&d-k=6,\\
t^{-\ff{d-k+2}{2(d-k-2)}},&d-k>6.
\end{cases}
$$
Moreover, letting     $\rr_\pp(x)$  be  the distance of $x$ to $\pp M$,     for large $t>0$   we have
$$
\W_2(\mu_t^{\dd_x},\mu_0)\asymp
\begin{cases}
t^{-1},&d\le 3,\\
t^{-1}\log (1+\rr_\pp(x)^{-1}), & d=4,\\
t^{-1} \rr_\pp(x)^{-1}, & d=5,\\
t^{-1}\sqrt{\log t} \rr_\pp(x)^{-1},&d=6,\\
t^{-\ff{d+2}{2(d-2)}}\rr_\pp(x)^{-\ff 4{d-2}},&d\ge 7.
\end{cases}
$$

\end{abstract}


 \vskip 2cm

\section{Introduction}

Let $M$ be a compact connected  $d$-dimensional Riemannian manifold with smooth boundary $\pp M$,   let $\rr$ be the Riemannian  distance on $M$,   let
$\scr P$ denote the set of all probability measures on $M$.  %For any $p\in [1,\infty)$,   the $p$-
The quadratic Wasserstein distance on $\scr P$ is defined as
$$
\W_2(\nu_1,\nu_2) :=  \inf_{\pi\in\C(\nu_1,\nu_2)} \left(\int_{M\times M}\rr(y,z)^2\,\pi(\d y,\d z)\right)^{\ff 1 2},\ \ \nu_1,\nu_2\in \scr P,
$$
where $\C(\nu_1,\nu_2)$ is the set of all  couplings for $\nu_1$ and $\nu_2$.

Let $V\in C^2(M)$ be such that $\mu(\d x)=\e^{V(x)}\d x\in \scr P$, where $\d x$ is the Riemannian volume measure.  Let $\DD$ and $\nn$ be the Laplacian and gradient operators on $M$.  We consider the     diffusion process $(X_t)_{t\ge0}$ generated by $L:=\DD+\nn V$
in $M^\circ:= M\setminus \pp M$ with life time
$$\tau:=\inf\big\{t\ge 0:X_t\in\pp M\big\}.$$
To ensure $\P(\tau>0)>0$, it is essential that    the initial distribution   belongs to the space
$$\scr P_0:= \big\{\nu\in \scr P:\ \nu(  M^\circ)>0\big\}.$$
In this case,   we consider the  empirical measure
$$\mu_t:=\ff1t\int_0^t\dd_{X_s}\d s,\ \ \ t\in (0,\tau),$$
where $\dd_{X_s}$ is the Dirac measure at point $X_s$,
and study the
long time behavior of the conditional  empirical measure
\beq\label{CM}
\mu_t^\nu:=\E^\nu\left(\mu_t \middle|\ t<\tau \right), \ \   \ \ \ t>0, \ \nu\in \scr P_0,
\end{equation}
 where   $\E^\nu$ is the expectation for the diffusion process $X_t$ with initial distribution $\nu$.

 According to \cite[Theorem 2.1(ii)]{CJ},  we have
 $$\lim_{t\to\infty} \|\mu_t^\nu-\mu_0\|_{var}=0,\ \ \ \nu\in \scr P_0,$$
where $\|\cdot\|_{var}$ is the total variation norm, and $$\mu_0(\d x):= \phi_0(x)^2\mu(\d x) $$ for
    the first unitary Dirichlet eigenfunction   $\phi_0 $ of $L$ in $L^2(\mu)$.  Without loss of generality, we assume that $\phi_0>0$ in $M^\circ$.

Let $\{ \ll_m\}_{m\ge 0}$ be all Dirichlet eigenvalues of $-L$ on $M$ listed in the increasing order including multiplicities, and let  $\{\phi_m\}_{m\ge 0}$ be the
associated unitary Dirichlet eigenfunctions, i.e.
$$
 L\phi_m=-\ll_m\phi_m,\ \ \  \phi_m|_{\pp M}=0,\ \mu(\phi_m^2)=1,\ m\ge 0.
$$
It is classical that,  see e.g., \cite{Chavel},   $\ll_0>0, \ll_1-\ll_0>0$ and
  $\{\phi_m\}_{m\ge 0}$ is an orthonormal basis of $L^2(\mu)$.
   Moreover, by \cite[Lemma 11(2)]{FS},   there exists a constant $\kk>1$ such that
 \beq\label{KK} \kk^{-1} \rr_\pp (x)\le \phi_0(x)\le \kk \rr_\pp(x),\ \ \ x\in M,\end{equation}
 where
 $$\rr_\pp (x):=\inf_{y\in\pp M} \rr(x,y)$$   is the  distance from $x$ to the boundary $\pp M$.

 According to \cite[Theorem 1.1]{W21}, for any $\nu\in \scr P_0(M)$, we have
\beq\label{0} \lim_{t\to\infty} \big\{t\W_2(\mu_t^\nu,\mu_0)\big\}^2=\ff 1 {\mu(\phi_0)^2\nu(\phi_0)^2} \sum_{m=1}^\infty\ff{\{\nu(\phi_0)\mu(\phi_m)+ \mu(\phi_0)\nu(\phi_m)\}^2}{(\ll_m-\ll_0)^3}>0.\end{equation}
 By  the Sobolev inequality, the limit is finite if either $d<6$, or $d\ge 6$  and $\rr_\nu:=\ff{\d\nu}{\d\mu}\in L^{p} (\mu) $ for some $p\in [\ff{2d}{d+6},\infty)\cap (1,\infty)$.
 The original   result claims that the limit is finite for  any $\nu\in \scr P_0$ if $d=6$, which is however not true since the Sobolev inequality (see, e.g., \cite[Theorem 2.3]{AB15})
 $$\big\|(-L)^{-\ff 3 2} f\big\|_{L^{\ff{dp}{d-3p}}(\mu)}\le K \|f\|_{L^p(\mu)},\ \ f\in L^p(\mu)$$
   %in \cite[(1.2)]{W21}
 only holds for $p\in (1,\ff d 3)$,  but not for $p=1$ which was used in the proof of  \cite[Theorem 1.1]{W21} for the case $d=6$. See also   \cite{WJEMS} for the sharp convergence rate of  $\E(\W_2(\mu_t,\mu_0)|\tau>t)$.

 There is also   a plentiful literature on the  Wasserstein convergence rate for empirical measures of   i.i.d. sequences and stochastic processes without killing,
see
   \cite{LW, LW24, TWZ, W25, W26, W, WW, WWZ, WZ} for  ergodic Markov processes, see \cite{HM, LW26} for   (subordinated) fractional Brownian motion on the flat tours, and see \cite{AT, BL, FG, HGT, LZ, WZ25} and references therein for the study of i.i.d. sequences.

However, when the  series in \eqref{0} diverges, it is the case when $d\ge 6$ and $\nu$ is singular, the sharp convergence rate of  $\W_2(\mu_t^\nu,\mu_0)$ remains  open.
In this paper, we aim to fill up this gap by presenting  the sharp convergence rate for any initial distribution $\nu\in\scr P_0$ in any  dimension $d$.

To state our main result, we first introduce some notation. Let
$$
L_0:=L+2\nn\log\phi_0=\Delta+\nn(V+2\log\phi_0),
$$ so that $(L_0, C_0^\infty(M^\circ))$ is symmetric in $L^2(\mu_0)$, whose closure $(L_0, \D(L_0))$ is a self-adjoint operator.
As shown in \cite{TW98} that the diffusion process generated by $L_0$ is
non-explosive in $M^\circ$.
Let $(Q_t)_{t\ge0}$ be the associated  symmetric diffusion semigroup, and let
  $q_t(x,y)$ be the  heat kernel of $Q_t$ with respect to $\mu_0$. Consider
 the Green function
\beq\label{GREEN}
G_0(x,y):=\int_0^\infty\bigl(q_t(x,y)-1\bigr)\d t,
\qquad x,y\in M^\circ.
\end{equation}
For any $\nu\in\scr P_0$, let
\beq\label{V0}
\nu_0(\d x):=\ff{\phi_0(x)}{\nu(\phi_0)}\nu(\d x),\qquad
g^\nu:=\int_{M^\circ}G_0(x,\cdot)\nu_0(\d x).
\end{equation}
We have  $g^\nu\in L^1(\mu_0)$
and  $\mu_0(g^\nu)=0$; see \eqref{G-REP} below.
Moreover, let
\beq\label{D-nu}
\mathcal D_\nu(t):=4\int_0^\infty s\,\mu_0\left(\ff{|\nn Q_sg^\nu|^2}{(1+t^{-1}|Q_sg^\nu|)^2}\right)\d s,\ \ \ t\ge 1.
\end{equation}
The main result of the paper is the following.

\beg{thm}\label{T1} There exist $C, t_0  \in (1,\infty)$ such that
for any $\nu\in\scr P_0$ and $t\ge t_0$, we have $\mathcal D_\nu(t)<\infty$  and
\beq\label{MAIN-D}
C^{-1}\ss{1+ \mathcal D_\nu(t) } \le  t \W_2(\mu_t^\nu,\mu_0) \le C \ss{1+ \mathcal D_\nu(t)},\ \ \ t\ge t_0,\  \nu\in\scr P_0.
\end{equation}
\end{thm}

To calculate $\mathcal D_\nu(t)$, we consider the following class of initial distributions for $k\in [0,d]$:
$$\scr P_{0,k}:=\big\{\nu\in \scr P_0: \
 \nu(B(x,r) ) \asymp r^k, \ \  x\in S_\nu:=M^\circ\cap {\rm supp}\nu,\ r\in [0,1]\big\},$$
where $B(x,r)$ is the geodesic ball in $M$ with center $x$ and  radius $r$, and   $A\asymp B$ for two positive variables $A$ and $B$  means that
$$ c^{-1}B\le A\le cB $$  holds for some constant $c\in (1,\infty)$.

It is clear that $\scr P_{0,0}$ contains $\nu\in \scr P_0$ with  finite support. When $k\in (0,d]$, $\nu\in \scr P_{0,k}$ means that it has   a local version of Ahlfors--David $k$-regularity, see \cite[Definition 1.1]{MS}.   In particular, for $k\in [1,d]\cap\mathbb N,$ $\scr P_{0,k}$ includes measures $\nu\in \scr P_0$ comparable with the volume measure of a $k$-dimensional submanifold  of $M$, see \cite[Proposition~1.12]{RZ}.

We have the following consequence of Theorem \ref{T1}, where constants $C_{\nu,n}$ and $t_\nu$ are  given in Proposition \ref{T2} below.

\beg{cor}\label{C1}  For any $k\in [0,d]$ and $\nu\in \scr P_{0,k}$, let $n=d-k$. Then there exist constants
 $ C_{\nu,n},  t_\nu\in   (1,\infty)$  such that for   any  $t\ge t_\nu$
 $$C_{\nu, n}^{-1} \vv_n(t)\le \W_2(\nu_t^\nu,\mu_0)\le C_{\nu,n} \vv_n(t),\ \ \ t\ge t_\nu$$ holds for
 $$\vv_n(t):=
\begin{cases}
t^{-1},&n<6,\\
 t^{-1}\ss{\log t},&n=6,\\
 t^{-\ff{n+2}{n-2}},&n>6.
\end{cases}$$
\end{cor}

Finally, we consider
  $\nu=\dd_x$, $x\in M^\circ$, and  denote  $\mu_t^x:=\mu_t^{\delta_x}$.
For any $\aa\in (6,\infty)$, let
$$t_{\aa, x}:=   \beg{cases} \rr_\pp(x)^{-d}, &d\notin\{4,6\},   \\
 \rr_\pp(x)^{-4 } \{\log(1+\rr_\pp(x)^{-1} )\}^{-2}, &d=4,\\
 \rr_\pp(x)^{-\aa}, &d=6,\end{cases}$$ which depends on $\aa$ only when $d=6$.
 For any $x\in M^\circ$ and $t>0$,  let
 $$\vv_{x,d}(t):= \begin{cases}
t^{-1},&d\le 3,\\
t^{-1}\log (1+\rr_\pp(x)^{-1}), & d=4,\\
t^{-1} \rr_\pp(x)^{-1}, & d=5,\\
t^{-1}\sqrt{\log t}\, \rr_\pp(x)^{-1},&d=6,\\
t^{-\ff{d+2}{2(d-2)}}\rr_\pp(x)^{-\ff 4{d-2}},&d\ge 7.
\end{cases}$$
We have the following result.

\beg{cor}\label{C2}  For any $\aa\in (1,\infty)$, there exists a constant $c\in (1,\infty)$ which depends on $\aa$ only when $d=6$, such that
$$\W_2(\mu_t^{x},\mu_0)\le c \vv_{x,d}(t),\ \ \ t\ge c,\ x\in M^\circ,$$
$$\W_2(\mu_t^{x},\mu_0)\ge c^{-1}  \vv_{x,d}(t),\ \ \ t\ge c\,t_{\aa, x} ,\ x\in M^\circ.$$

\end{cor}

\section{Some preparations}

Let  $p_t^D(x,y)$ be the Dirichlet heat kernel of $L:=\DD+\nn V$   with respect to $\mu$,  so that the  associated  Dirichlet semigroup $P_t^D$    is formulated as
$$
P_t^Df(x) :=\E^x\left[f(X_t)1_{\{t<\tau\}}\right]=\int_Mp_t^D(x,y)f(y)\,\mu(\d y),\ \ x\in M,\ t\ge 0,\ f\in \B_b(M),
$$ where $\B_b(M)$ is the set of all bounded measurable functions on $M$.
By the ground state transform, we have
\beq\label{QT}
Q_tf=\e^{\ll_0t}\phi_0^{-1}P_t^D(f\phi_0),
\qquad f\in L^2(\mu_0),\quad t>0,
\end{equation}
so that
\beq\label{QK}
q_t(x,y)=
\ff{\e^{\ll_0t}p_t^D(x,y)}{\phi_0(x)\phi_0(y)},
\qquad x,y\in M^\circ,\quad t>0.
\end{equation}
\black
Moreover,
$ \{\psi_m:=\phi_m \phi_0^{-1}\}_{m\ge 0} $ is an orthonormal basis of  $L^2(\mu_0)$ such that
\beq\label{L0} L_0 \psi_m= -(\ll_m-\ll_0)\psi_m,\ \ \ m\ge 0.\end{equation}
By the spectral representation, we have
\beq\label{QE}\beg{split}
&Q_t \psi_m =\e^{-(\ll_m-\ll_0)t}\psi_m,\qquad\ t\ge 0,\ m\ge 1,\\
&\|Q_t-\mu_0 \|_{L^2(\mu_0)}:=\sup_{\|f\|_{L^2(\mu_0)}\le 1}\|Q_tf-\mu_0 (f)\|_{L^2(\mu_0)} \le \e^{-(\ll_1-\ll_0)t},\ \ \ t\ge 0,\\
&q_t(x,y)=1+ \sum_{m=1}^{\infty}\e^{-(\ll_m-\ll_0)t}\psi_m(x)\psi_m(y),\ \ \ t>0,\ x,y\in M^\circ.\end{split}
\end{equation}

We will need  the following  heat kernel estimate  which is due to  \cite{Zhang} when $V=0$.

\beg{lem}\label{HK} There exists a constant $k>1$ such that
\beq\label{DU} \beg{split}& \bigg(\ff{\rr_\pp(x)}{(1\land  t)^{1/2}}\land 1\bigg) \bigg(\ff{\rr_\pp(y)}{ (1\land t)^{1/2}}\land 1\bigg)
\ff{1} {k (1\land   t)^{ d/2}} \e^{-\ll_0t - k\rr(x,y)^2/t}\le p_t^D(x,y)\\
&\quad \le \bigg(\ff{\rr_\pp(x)}{(1\land  t)^{1/2}}\land 1\bigg) \bigg(\ff{\rr_\pp(y)}{ (1\land t)^{1/2}}\land 1\bigg)
\ff{k} {(1\land  t)^{ d/2}} \e^{ -\ll_0t - \rr(x,y)^2/(kt)},\ \ \ t>0,\ x,y\in M^\circ.\end{split}\end{equation}
Consequently, there exists a constant  $k'>1$ such that
\beq\label{UC}  \beg{split}  |q_t (x,y)-1|      \le k' \min\bigg\{  \ff{(1\land t)^{-\ff d 2}}{ \rr_\pp(x) \rr_\pp(y)},\   (1\land t)^{-\ff{d+2}2}  \e^{- t (\ll_1-\ll_0)}\bigg\},\ \
  t>0,\ x,y\in M^\circ,\end{split} \end{equation}
\beq \label{DHL}
q_t(x,x)\ge  \ff {t^{-d/2}}  {k'(t+ \rr_\pp(x)^2)}     \ge \ff {t^{-d/2}} {k'(1+\|\rr_\pp\|_\infty^2)},\qquad  \ \ t\in (0,1],\ x\in M^\circ.
\end{equation}
 \end{lem}

\beg{proof} When $V=0$, \eqref{DU} is included in \cite[Theorem 1.1(b)]{Zhang}. It is easy to see that the same estimate holds for $V\in C^2(M).$ Below we give a proof by using
the unitary transform
$$
\mathcal U: L^2(M,\mu)\to L^2(M,\d x), \qquad \mathcal Uf:=\e^{V/2}f.
$$

(a) We first prove \eqref{DU} for short time, i.e., for any $t_0\in (0,\infty)$
 there exists a constant $k_0>1$ such that
\beq\label{DUU} \beg{split}& \Big(\ff{\rr_\pp(x)}{  t^{1/2}}\land 1\Big) \Big(\ff{\rr_\pp(y)}{  t^{1/2}}\land 1\Big)
\ff{1} {k_0   t^{ d/2}} \e^{- k_0\rr(x,y)^2/t}\le p_t^D(x,y)\\
&\quad \le \Big(\ff{\rr_\pp(x)}{  t^{1/2}}\land 1\Big) \Big(\ff{\rr_\pp(y)}{  t^{1/2}}\land 1\Big)
\ff{k_0} {  t^{ d/2}} \e^{ -\ll_0t - \rr(x,y)^2/(k_0t)},\ \ \ t\in (0,t_0],  \ x,y\in M^\circ.\end{split}\end{equation}
By a direct calculation, we have
\beq\label{CONJ}
\mathcal U L \mathcal U^{-1}=\DD -W\ \text{on}\ C_b^2(M),\qquad W:=\ff12\DD V+\ff14|\nn V|^2\in C_b(M).
\end{equation}
  Let $ X_t^0$ be the   diffusion process generated by $\DD$ with lifetime
 $$\tau^0:= \inf\big\{t\ge 0: X_t^0\in\pp M\big\}.$$
By  the Feynman--Kac formula, for any  $f\in \B_b^+(M):=\{f\in \B_b(M),\ f\ge 0\}$,    we have
$$ \e^{V(x)/2}P_t^D(\e^{-V/2}f)(x)=\E^x\left[1_{\{t<\tau^0\}} f(X_t^0)\e^{-\int_0^t W(X_s^0)\d s} \right],\quad t>0.$$
So, letting $p_t^0(x,y)$ be the  Dirichlet
heat kernel of  $\DD $
with respect to the volume measure,   we obtain
\beg{align*}&\e^{-\|W\|_\infty t} \int_M p_t^0(x,y) f(y) \d y \le  \e^{V(x)/2} \int_M p_t^D(x,y) \e^{V(y)/2} f(y) \d y \\
&\qquad  \le \e^{\|W\|_\infty t} \int_M p_t^0(x,y) f(y) \d y,\ \ \  f\in \B_b^+(M).\end{align*}
This implies
 $$\e^{-(\|W\|_\infty +  \|V\|_\infty) t_0}   p_t^0(x,y)  \le  p_t^D(x,y) \le p_t^0(x,y)\e^{(\|W\|_\infty  + \|V\|_\infty)t_0},\ \ \ t\in (0, t_0].  $$
By \cite[Theorem 1.1(b)]{Zhang},  $p_t^0(x,y)$ satisfies \eqref{DUU} for some constant $t_0>0$, so $p_t^D(x,y)$ satisfies the same type estimate with a different constant $k_0>1$.

(b)  By \eqref{KK}, \eqref{QK} and \eqref{DUU},   \eqref{UC}  and  \eqref{DHL} hold  for  $t\in (0,1]$.  We intend  to prove
 \eqref{UC}   for $t>1$, i.e. there exists a  constant $k'>1$ such that
 \beq\label{UC'}  |q_t(x,y)-1|\le  k' (1\land t)^{-\ff{d+2}2}  \e^{- t (\ll_1-\ll_0)},\ \ \ t>1,\ x,y\in M^\circ.\end{equation}
  By \eqref{UC} for $t=1/2$, we have
 $$\|Q_{1/2}-\mu_0\|_{L^2(\mu_0)\to L^\infty(\mu_0)}\le \|Q_{1/2}-\mu_0\|_{L^1(\mu_0)\to L^\infty(\mu_0)}\le k'2^{\ff{d+2}2},$$
 which together with  the second inequality in \eqref{QE}  yields
\beg{align*}&|q_t(x,y)-1|= \big|(Q_{1/2}-\mu_0)(q_{t-1/2}(x,\cdot) -1)(y) \big|\le k' 2^{\ff{d+2}2} \|q_{t-1/2}(x,\cdot) -1\|_{L^2(\mu_0)}\\
&= k' 2^{\ff{d+2}2} \big\|(Q_{t-1}-\mu_0) (q_{1/2}(x,\cdot) -1)\big\|_{L^2(\mu_0)}\le k' 2^{\ff{d+2}2} \e^{-(t-1)(\ll_1-\ll_0)} \ss{q_1(x,x)-1}.\end{align*}
Combining this with \eqref{UC} for $t=1$, we prove \eqref{UC'} for   a different constant $k'>1$.

(c) It remains to prove \eqref{DU} for large time. By   \eqref{UC}, there exists a constant $t_0\ge 1$ such that
$$\|q_t -1\|_\infty\le  \ff 1 2,\ \ \ t\ge t_0,$$
so that  \eqref{QK} and \eqref{KK}  imply
$$\ff 1 2 \e^{-\ll_0 t} \rr_\pp(x)\rr_\pp(y)  \le p_t^D(x,y)\le \ff 3 2 \e^{-\ll_0t} \rr_\pp(x)\rr_\pp(y),\ \ \ \ t\ge t_0,\ x,y\in M.$$
Combining this with
$$\ff{\rr_\pp(x)}{1\lor \|\rr_\pp \|_\infty}\le \ff{\rr_\pp(x)}{(1\land t)^{1/2}}\land 1 \le \rr_\pp(x),\ \ \ \ t\ge t_0\ge 1,$$ we derive
\eqref{DU} for  some constant $k>1$ and all $t\ge t_0$.  \end{proof}

The following result is due to   \cite[Theorem A.1]{W22},  which extends  \cite[Theorem 2]{Ledoux}.

\beg{lem}
Let $\d\nu_i=f_i\d\mu_0\in\scr P$, $i=1,2$.
If $f_1+f_2\ge\theta$ for some constant $\theta>0$, then
\beq\label{WS}
\W_2(\nu_1,\nu_2)^2
\le \min\left\{8\theta^{-1}\|\nu_1-\nu_2\|_{H^{-1}(\mu_0)}^2,4\mu_0\left(
\ff{|\nn(-L_0)^{-1}(f_1-f_2)|^2}{f_1}\right)\right\}.
\end{equation}
\end{lem}

To apply \eqref{WS}, we need to calculate the density function of $\mu_t^\nu$ with respect to $\mu_0$. To this end,
we introduce
    \beq\label{WST} \mu_{s,t}^\nu:= \E^\nu(\dd_{X_s}\mid t<\tau),\ \ \ \nu\in \scr P_0,\ 0\le s\le t<\infty,\end{equation}
so that
\beq\label{MUT} \mu_t^\nu =\ff 1 t \int_0^t \mu_{s,t}^\nu\d s,\ \ \ \nu\in\scr P_0,\ t>0.\end{equation}
With \eqref{V0} let
\beq\label{QS} q_s^\nu(y):=\nu_0(q_s(\cdot,y))= \ff 1 {\nu(\phi_0)} \int_{M}  \phi_0 q_s(\cdot,y)\d\nu,\ \ \  y\in M^\circ.\end{equation}

\beg{lem}
For any   $\nu\in\scr P_0$,
\beq\label{CDF}
h_{s,t}^\nu(y):=\ff{\d\mu_{s,t}^\nu }{\d\mu_0}(y)=\ff{q_s^\nu(y)Q_{t-s}\phi_0^{-1}(y)}{\nu_0(Q_t\phi_0^{-1})},\qquad y\in M^\circ,\ 0<s\le t<\infty.
\end{equation}
  Moreover, there exists a constant $s_0\in [2,\infty)$ such that
\beq\label{QB}
\ff 12 \mu_0(\phi_0^{-1}) \le \nu_0(Q_t\phi_0^{-1})\le  2 \mu_0(\phi_0^{-1}), \qquad t\ge s_0,\ \nu\in \scr P_0.
\end{equation}
\end{lem}

\beg{proof}

 By  the Markov property of $X_t$, we have
$$
  \E^\nu\left[f(X_s)1_{\{t<\tau\}}\right] = \E^\nu\left[f(X_s)1_{\{s<\tau\}}\E^{X_s}  1_{\{t-s<\tau\}}  \right]=\nu(P_s^D\left(fP_{t-s}^D1\right)),\ \ f\in \B_b (M).
$$
Combining this with  \eqref{QT} and \eqref{V0},   we obtain
\beg{align*}& \mu_{s,t}^\nu(f):=    \ff{\E^\nu\left[f(X_s)1_{\{t<\tau\}}\right]}{\P^\nu(t<\tau)}=  \ff{\nu(P_s^D\left(fP_{t-s}^D1\right))}{\nu(P_t^D1)}\\
 &= \ff{\int_M\e^{-\ll_0t}\phi_0(x) Q_s\left(fQ_{t-s}\phi_0^{-1}\right)(x)\nu(\d x)}{\int_M \e^{-\ll_0t}\phi_0(x)Q_t\phi_0^{-1}(x)\nu(\d x)}\\
 &= \ff{\nu_0\left(  Q_s\left(fQ_{t-s}\phi_0^{-1}\right)\right)}{ \nu_0( Q_t\phi_0^{-1})},\ \ \ f\in \B_b(M).\end{align*}
Noting that $q_s(x,y)$ is the heat kernel of $Q_s$ with respect to $\mu_0$,  this together with \eqref{QS} implies   \eqref{CDF}.

 Next, by $\mu_0(\phi_0^{-2})=\mu(M)=1$ and Jensen's inequality, we have $\mu_0(\phi_0^{-1})\in (0,1]$. Then    \eqref{UC}  implies
 $$\|Q_t\phi^{-1}_0 -\mu_0(\phi_0^{-1})\|_{\infty}\le c \e^{- t(\ll_1-\ll_0)}\mu_0(\phi_0^{-1}),\ \ \ t\ge1.$$
   Taking $s_0\ge2$ sufficiently large, we obtain
$$
\ff 12 \mu_0(\phi_0^{-1}) \le Q_t\phi_0^{-1}(x)\le  2 \mu_0(\phi_0^{-1}), \qquad t\ge s_0,\quad x\in M^\circ.
$$
Noting that $\nu_0$ is a probability measure, we  obtain  \eqref{QB}.
\end{proof}

By \eqref{UC}  and $\|\psi_m\|_{L^2(\mu_0)}=1,$ we find a constant $k_1\in (0,\infty)$ such that
 $$\|Q_t \psi_m\|_\infty \le  \bigg(\int_M q_{2t}(x,x) \mu_0(\d x)\bigg)^{\ff 1 2} \le k_1 (1\land t)^{-\ff{d+2}4},\ \ t>0,\ m\ge 0.$$
Combining this with  \eqref{QE}, we find a constant $k_2\in (0,\infty)$ such that
\beq\label{QMN} \|\psi_m\|_\infty=\inf_{t>0} \e^{(\ll_m-\ll_0)t} \|Q_t \psi_m\|_\infty\le k_2 (\ll_m-\ll_0)^{\ff{d+2}4},\ \ m\ge 1.\end{equation}
  For any finite signed measure   $\zeta$ on $M$, let $|\zeta|=\zeta^+-\zeta^-$ be its total variation, and denote
$$\zeta(f)= \int_{M} f \d\zeta^+- \int_M f\d\zeta^-,\ \ \ f\in L^1(|\zeta|).$$
By \eqref{QMN}, we have $|\zeta(\psi_m)|\le |\zeta|(M) \|\psi_m\|_\infty<\infty$ for any finite signed measure $\zeta$ on $M$.
If $\zeta(M)=0$,  then   \eqref{L0} implies
\beq\label{HN} \|\zeta\|_{H^{-1}(\mu_0)} := \sup_{\mu_0(f^2)\le 1} \zeta\Big((-L_0)^{-\ff 1 2} f\Big)  = \bigg(\sum_{m=1}^{\infty}\ff{|\zeta(\psi_m)|^2}
{\ll_m-\ll_0}\bigg)^{\ff 1 2}\in [0,\infty].
\end{equation}
 When $\ff{\d\zeta}{\d\mu_0}$ exists, we also denote
$$\|\zeta\|_{H^{-1}(\mu_0)}= \bigg\|\ff{\d\zeta}{\d\mu_0}\bigg\|_{H^{-1}(\mu_0)}.$$

\beg{lem}\label{DISP}
There exists a constant $c_\rr>0$ such that
\beq\label{QM}
Q_s\rr(x,\cdot)^2(x)\le c_\rr (1\land s),\qquad x\in M^\circ,\quad s>0.
\end{equation}
Moreover, for any $\nu\in \scr P_0$,
\beq\label{D-INF}
\mathcal D_\nu(\infty) :=\lim_{t\to\infty}\mathcal D_\nu(t)=4\int_0^\infty s\,\mu_0\bigl(|\nn Q_sg^\nu|^2\bigr)\d s
 =\|g^\nu\|_{H^{-1}(\mu_0)}^2.
\end{equation}
\end{lem}

\beg{proof}  By \eqref{KK}, \eqref{QT} and \eqref{DU}, we find constants $c_1,c_2>0$ such that
\beq\label{QS1} \beg{split}
&Q_s\rr(x,\cdot)^2(x) =\ff{\e^{\ll_0s}}{\phi_0(x)}\int_M\rr(x,y)^2\phi_0(y)p_s^D(x,y)\mu(\d y)\\
&\le \ff{c_1 s^{-d/2}}{\phi_0(x)}\left(1\wedge\ff{\phi_0(x)}{\sqrt{s}}\right)\int_M \rr(x,y)^2(\phi_0(x)+\|\nn \phi_0\|_\infty\rr(x,y))\e^{-\rr(x,y)^2/(c_1s)}\mu(\d y)\\
&\le\ff{c_2 s^{-d/2}}{\phi_0(x)+\sqrt s}
\int_M\rr(x,y)^2\big(\phi_0(x)+\rr(x,y)\big)\e^{-\rr(x,y)^2/(c_1s)}\mu(\d y),\ \ \ s\in (0,1].
\end{split}\end{equation}
Noting that  $\mu(B(x,r))\le Cr^d$ holds for some constant $C\in (0,\infty)$ and all $r>0$, we find  constant $c_3>0$ such that for  $j=2,3$,
\beg{align*}
&\int_M\rr(x,y)^j\e^{-\rr(x,y)^2/(c_1s)}\mu(\d y) =\sum_{m=0}^\infty \int_{\{m\sqrt{s}\le\rr(x,y)\le(m+1)\sqrt{s}\}}\rr(x,y)^j \e^{-\rr(x,y)^2/(c_1 s)}\mu(\d y)\\
&\le s^{j/2}\sum_{m=0}^\infty(m+1)^j\e^{-m^2/c_1}\mu\big(B(x,(m+1)\sqrt s)\big)
 \le c_3s^{(d+j)/2},\qquad s\in (0,1].
\end{align*}
Combining this with \eqref{QS1}, we find   a constant $c_4 \in (0,\infty)$ such that
$$
Q_s\rr(x,\cdot)^2(x)\le\ff{c_4(\phi_0(x)s+s^{3/2})}{\phi_0(x)+\sqrt s}=c_4s,\qquad s\in (0,1],\ x\in M^\circ.
$$ Since $\rr$ is bounded,  this implies   \eqref{QM}  for some constant $c_\rr$.

 Next, since  $q_t$ is the heat kernel of the Markov semigroup $Q_t$ with respect to $\mu_0$, we have   $q_t\ge 0$ and  $\mu_0(q_t(y,\cdot))=1$. So,
  \eqref{UC}  implies
$$
c:=\sup_{y\in M^\circ}\int_0^\infty\mu_0\bigl(|q_t(y,\cdot)-1|\bigr)\d t\le 2+\int_1^\infty\sup_{x,y\in M^\circ}|q_u(x,y)-1|\d u <\infty,$$
and the Green potential $g^\nu$ defined in \eqref{V0} satisfies
\beq\label{G-REP}
 \mu_0(|g^\nu|) \le c,\quad \mu_0(g^\nu)=0,\qquad \nu\in \scr P_0.
\end{equation}
This together with \eqref{QMN} implies
\beq\label{QM*} \mu_0(g^\nu \psi_m)^2\le (ck_2)^2 (\ll_m-\ll_0)^{\ff{d+2}2},\ \ \ m\ge 0.\end{equation}
Moreover, by   \eqref{UC} and \eqref{G-REP}  we have  $Q_tg^\nu\in L^\infty(\mu_0)\subset L^2(\mu_0)$ for $t>0$. Combining this with   $\eqref{QE}$ and the symmetry of $Q_t$ in $L^2(\mu_0)$, we have
%, so that the expansion with orthonormal basis $\{\psi_m\}_{m\ge 0}$ gives
\beq\label{QS**} Q_t g^\nu= \sum_{m=0}^\infty \mu_0(\psi_m Q_tg^\nu) \psi_m=\sum_{m=0}^\infty \mu_0(g^\nu\psi_m) \e^{-(\ll_m-\ll_0)t} \psi_m.\end{equation}
%This together with   $\eqref{QE}$ and the symmetry of $Q_t$ in $L^2(\mu_0)$  gives
%$$Q_t g^\nu= \sum_{m=0}^\infty \mu_0(g^\nu\psi_m) \e^{-(\ll_m-\ll_0)t} \psi_m.$$
Noting that
  $$\mu_0(\<\nn \psi_m,\nn \psi_n\>)    =  1_{\{n=m\}} ( \ll_m-\ll_0) \asymp  1_{\{n=m\}} m^{\ff 2 d},\ \ \ m\ge 1,$$
so that \eqref{QM*}   ensures
$$\sum_{m=0}^\infty  \mu_0(g^\nu\psi_m)^2 \e^{-2(\ll_m-\ll_0)t} \mu_0(|\nn \psi_m|^2)   <\infty,\ \ \ t>0,$$
we deduce from \eqref{QS**} that $Q_tg^\nu\in \D((-L_0)^{\ff 1 2})$ and
$$\mu_0\bigl(|\nn Q_tg^\nu|^2\bigr) =\| (-L_0)^{\ff 1 2} Q_t g^\nu\|_{L^2(\mu_0)}^2 =\sum_{m=1}^\infty(\ll_m-\ll_0)\e^{-2t(\ll_m-\ll_0)}\mu_0(g^\nu\psi_m)^2<\infty,\qquad t>0.
$$ Therefore,
\begin{align*}
&\mathcal D_\nu(\infty):= 4\int_0^\infty s\,\mu_0\bigl(|\nn Q_sg^\nu|^2\bigr)\d s\\
&=4\sum_{m=1}^\infty(\ll_m-\ll_0)\mu_0(g^\nu\psi_m)^2\int_0^\infty s\e^{-2(\ll_m-\ll_0)s}\d s =\sum_{m=1}^\infty\ff{\mu_0(g^\nu\psi_m)^2}{\ll_m-\ll_0}=\|g^\nu\|_{H^{-1}(\mu_0)}^2.
\end{align*}
 \end{proof}

\section{Proof of  Theorem \ref{T1}: the upper bound estimate}

We first explain the main idea of the proof.

Firstly, by the exponential ergodicity of $Q_t$ we have $Q_t \phi_0^{-1}\to\mu_0(\phi_0^{-1})=\mu(\phi_0)\in (0,\infty)$ exponentially fast, so that
by \eqref{CDF}
$$h_{s,t}^\nu:=\ff{\d\mu_{s,t}^\nu}{\d\mu_0}\to q_s^\nu\  \text{as} \ t\to\infty,\ \text{exponentially\  fast}.$$    This and \eqref{MUT} suggest us to simplify $\mu_t^\nu$ as
\beq\label{MT-nu} \d\bar \mu_t^\nu:= g_t^\nu\d\mu_0,\ \ \ g_t^\nu:=  \ff 1 t \int_0^t q_s^\nu\d s,\ \ \ t>0.\end{equation}
Then, by the triangle inequality of $\W_2$,
\beq\label{UP1} \W_2(\mu_t^\nu,\mu_0) \le  \W_2(\mu_t^\nu,\bar \mu_t^\nu) +  \W_2(\bar\mu_t^\nu,\mu_0),\ \ \ t>0.\end{equation}
We first show that   $\W_2(\mu_t^\nu,\bar \mu_t^\nu)\le c t^{-1}$ for some constant $c\in (0,\infty)$,   then estimate $ \W_2(\bar\mu_t^\nu,\mu_0)$ by using \eqref{WS}.

\beg{lem}\label{JL}
There exists a  constant  $c\in (0,\infty)$  such that
\beq\label{QCOMP}
\W_2(\mu_t^\nu,\bar\mu_t^\nu) \le ct^{-1},\qquad t>0,\quad \nu\in\scr P_0.
\end{equation}
 \end{lem}

\beg{proof} Since $M$ is compact, there exists a constant $c_1\in (0,\infty)$ such that $\W_2(\mu_t^\nu,\bar\mu_t^\nu)\le c_1$, so it suffices
to find constants $c, t_0\in (0,\infty)$ such that \eqref{QCOMP} holds for any $\nu\in \scr P_0$ and  $t\ge t_0.$

By \eqref{UC}, there exists a constant $s_1>0$ such that
$$
q_s(x,y)\ge\ff12,\qquad s\ge s_1,\quad x,y\in M^\circ,
$$ so that \eqref{QS} and \eqref{MT-nu} imply
\beq\label{MD2-nu}
g_t^\nu\ge\ff{t-s_1}{2t}\ge\ff14,\qquad t\ge2s_1.
\end{equation}
Let $s_0\in [2,\infty)$ be in \eqref{QB} and define
\beq\label{T00} t_0:=\max\{2s_1,s_0+1\}.\end{equation}
Let
$$
\xi_t^\nu:=\ff1t\bar\mu_1^\nu+\ff1t\int_1^t\mu_{s,t}^\nu\d s,\ \ t>0.
$$
We have
\beq\label{UP2} \W_2(\mu_t^\nu, \bar\mu_t^\nu)\le  \W_2(\mu_t^\nu, \xi_t^\nu)+ \W_2(\xi_t^\nu,\bar\mu_t^\nu).\end{equation}
Below we estimate the two terms in the upper bound respectively.

Firstly,  since $D:=\|\rr\|_\infty<\infty,$ we obtain
\beq\label{UP3} \W_2(\mu_t^\nu, \xi_t^\nu)^2\le D^2 |\mu_t^\nu- \xi_t^\nu|(M)= \ff{D^2} t \int_0^1 \|h_{s,t}^\nu- q_s^\nu\|_{L^1(\mu_0)} \d s.\end{equation}
By \eqref{CDF},
\beq\label{*00}
h_{s,t}^\nu-q_s^\nu=\ff{q_s^\nu}{\nu_0(Q_t\phi_0^{-1})}\big(Q_{t-s}\phi_0^{-1}-\nu_0(Q_t\phi_0^{-1})\big).\end{equation}
 By \eqref{UC},
\beq\label{Qnu}
\|Q_s\phi_0^{-1}-\mu_0(\phi_0^{-1})\|_\infty\le c\e^{-s(\ll_1-\ll_0)},\qquad s\ge1.
\end{equation}
Since $\nu_0$ is a probability measure, this  implies
\begin{equation}\label{Qt}
\left|\nu_0(Q_t\phi_0^{-1})-\mu_0(\phi_0^{-1})\right|\le c\e^{-t(\ll_1-\ll_0)},\qquad t\ge1.
\end{equation}
 Combining this with \eqref{QB}, \eqref{Qnu}  and that  $q_s^\nu$ is a probability density with respect to $\mu_0$, we find a constant $c_1\in (0,\infty)$ such that
 \beg{align*}
\|h_{s,t}^\nu-q_s^\nu\|_{L^1(\mu_0)}&\le\ff{\|q_s^\nu\|_{L^1(\mu_0)}}{\nu_0(Q_t\phi_0^{-1})}
\left[\|Q_{t-s}\phi_0^{-1}-\mu_0(\phi_0^{-1})\|_\infty+\left|\nu_0(Q_t\phi_0^{-1})-\mu_0(\phi_0^{-1})\right|
\right]\\
&\le c_1\e^{-t(\ll_1-\ll_0)},\ \ \ s\in (0,1],\ t\ge t_0.\end{align*}
This together with \eqref{UP3} yields that for some constant $c_2\in (0,\infty)$
\beq\label{UP4} \W_2(\mu_t^\nu, \xi_t^\nu) \le c_2 t^{-1},\ \ \ t\ge t_0.\end{equation}

Next, since $\mu_0(\phi_0^{-2})=1$,   \eqref{QE} implies
\beq\label{Qnu2}
\|Q_s\phi_0^{-1}-\mu_0(\phi_0^{-1})\|_{L^2(\mu_0)}\le \e^{-s(\ll_1-\ll_0)},\qquad s\ge0.
\end{equation}
Combining this with \eqref{*00} and \eqref{Qt}, we find a constant $c_3\in (0,\infty)$ such that
\beg{align*}
\|h_{s,t}^\nu-q_s^\nu\|_{L^2(\mu_0)}
&\le\ff{\|q_s^\nu\|_\infty}{\nu_0(Q_t\phi_0^{-1})}\left[\|Q_{t-s}\phi_0^{-1}-\mu_0(\phi_0^{-1})\|_{L^2(\mu_0)}
+\left|\nu_0(Q_t\phi_0^{-1})-\mu_0(\phi_0^{-1})\right|\right]\\
&\le c_3\e^{-(t-s)(\ll_1-\ll_0)},\ \ \ t\ge s\ge 1.
\end{align*}
This together with \eqref{WS}, \eqref{MD2-nu} and the Poincar\'e inequality
$$\|f\|_{H^{-1}(\mu_0)}^2\le \ff 1 {\ll_1-\ll_0} \|f\|_{L^2(\mu_0)},\ \ f\in L^2(\mu_0),\ \mu_0(f)=0$$
leads to that for some $c_4\in (0,\infty)$
\beg{align*} &\W_2(\xi_t^\nu,\bar\mu_t^\nu)^2 \le 32 \|\xi_t^\nu-\bar\mu_t^\nu\|_{H^{-1}(\mu_0)}^2\\
&\le
\ff{32}{(\ll_1-\ll_0)t^2} \bigg(\int_1^t  \|h_{s,t}^\nu-q_s^\nu\|_{L^2(\mu_0)}\d s\bigg)^2\le c_4 t^{-2},\ \ t\ge t_0.\end{align*}
Then the proof is finished by combining this with \eqref{UP2} and \eqref{UP4}.
  \end{proof}

To estimate $\W_2(\bar\mu_t^\nu,\mu_0)$ using \eqref{WS}, we will need the following lemma.

\beg{lem}\label{FISH}
There exists a constant $c>0$ such that
\beq\label{FI}
\int_0^r s^2\mu_0\left(\ff{|\nn q_s(x,\cdot)|^2}{q_s(x,\cdot)}\right)\d s\le cr^2,\qquad x\in M^\circ,\quad r\in(0,1].
\end{equation}
\end{lem}

\beg{proof}
For $x\in M^\circ$ and $s>0$, let
$$
H_x(s):=\mu_0\big(q_s(x,\cdot)\log q_s(x,\cdot)\big).
$$
By \eqref{KK},
$$
\ff1{\kappa(\phi_0(x)+\sqrt s)}\le\ff{1\wedge(\rr_\pp(x)/\sqrt s)}{\phi_0(x)}\le\ff{2\kappa}{\phi_0(x)+\sqrt s},\qquad x\in M^\circ.
$$
Combining  this with
\eqref{QK} and \eqref{DU}, we find  a constant $c_1\in (1,\infty)$ such that
\beq\label{QGA}\beg{split}
&\ff{s^{-d/2}\e^{-c_1\rr(x,y)^2/s}}{c_1(\phi_0(x)+\sqrt s)(\phi_0(y)+\sqrt s)}\\
 &\le q_s(x,y)\le\ff{c_1s^{-d/2}\e^{-\rr(x,y)^2/(c_1s)}}{(\phi_0(x)+\sqrt s)(\phi_0(y)+\sqrt s)},\quad \ 0<s\le1,\ x,y\in M^\circ.\end{split}
\end{equation}
In particular,
$$
q_s(x,y)\le\ff{c_1s^{-d/2}}{(\phi_0(x)+\sqrt s)(\phi_0(y)+\sqrt s)}\le c_1s^{-(d+2)/2},\qquad 0<s\le1.
$$
Since $\mu_0(q_s(x,\cdot))=1$,  this implies
$$
0\le H_x(s)\le\log\|q_s(x,\cdot)\|_\infty\le\log c_1+\ff{d+2}2\log\ff1s,\qquad 0<s\le1.
$$
Consequently,
\beq\label{HZ}
\lim_{s\downarrow0}s^2H_x(s)=0.
\end{equation}
Noting that
$$\partial_s q_s(x,\cdot)=L_0q_s(x,\cdot),\ \ \
 \mu_0(L_0q_s(x,\cdot))=0,$$  using the integration by parts formula we obtain
\beg{align*}
&\ff{\d}{\d s}H_x(s) =\mu_0\big((1+\log q_s(x,\cdot))L_0q_s(x,\cdot)\big)
 =\mu_0\big(\log q_s(x,\cdot)L_0q_s(x,\cdot)\big)\\
&=-\mu_0\big(\langle\nn\log q_s(x,\cdot),\nn q_s(x,\cdot)\rangle\big)
 =-\mu_0\left(\ff{|\nn q_s(x,\cdot)|^2}{q_s(x,\cdot)}\right),\qquad s>0.
\end{align*}
This together with \eqref{HZ} yields
\beq\begin{split}\label{IU}
&\int_0^r s^2\mu_0\left(\ff{|\nn q_s(x,\cdot)|^2}{q_s(x,\cdot)}\right)\d s =-\lim_{\delta\downarrow0}\int_\delta^r s^2H_x'(s)\d s\\
&=-r^2H_x(r)+2\int_0^r sH_x(s)\d s
 =2\int_0^r s\big(H_x(s)-H_x(r)\big)\d s,\ \ \ r>0.\end{split}\end{equation}
To bound $H_x(s)-H_x(r)$, we estimate $H_x$ from above and below respectively.

Firstly, by Jensen's inequality,  \eqref{QGA} and noting that $\mu_0(q_s(x,\cdot)^2)=q_{2s}(x,x),$ we obtain
\beq\label{HS}
\beg{split}
H_x(s)&\le\log\mu_0\big(q_s(x,\cdot)^2\big)=\log q_{2s}(x,x)\\
&\le\log c_1-\ff d2\log s-2\log(\phi_0(x)+\sqrt s),\ \ \ s\in (0,1].
\end{split}
\end{equation}

Next, by the Lipschitz continuity of $\phi_0$ and \eqref{QM}, we obtain
\beg{align*}
Q_r\phi_0(x)&\le\phi_0(x)+\|\nn\phi_0\|_\infty\big(Q_r(\rr(x,\cdot)^2)(x)\big)^{1/2}\\
&\le\phi_0(x)+\|\nn\phi_0\|_\infty\sqrt{c_\rr r},\ \ \ r>0,\ x\in M^\circ.
\end{align*}
So, by Jensen's inequality, for $r>0$ and $x\in M^\circ$
\beg{align*}
&Q_r\big(\log(\phi_0+\sqrt r)\big)(x) \le\log\big(Q_r\phi_0(x)+\sqrt r\big)\le\log(\phi_0(x)+\sqrt r)+\log\big(1+\|\nn\phi_0\|_\infty\sqrt{c_\rr}\big).
\end{align*}
Taking  logarithm for  the lower bound in \eqref{QGA} and integrating against $q_r(x,y)\mu_0(\d y)$,  and using \eqref{QM} again, we find a constant $c_2\in (0,\infty)$ such that
\beg{align*}
H_x(r)&\ge-\ff d2\log r-\log(\phi_0(x)+\sqrt r)-Q_r\big(\log(\phi_0+\sqrt r)\big)(x)\\
&\quad-\log c_1-\ff{c_1}r Q_r(\rr(x,\cdot)^2)(x)\\
&\ge-c_2-\ff d2\log r-2\log(\phi_0(x)+\sqrt r),\ \ r\in (0,1],\ x\in M^\circ
\end{align*}
Combining this with
$$
\ff{\phi_0(x)+\sqrt r}{\phi_0(x)+\sqrt s}=\sqrt{\ff rs}\,\ff{1+\phi_0(x)/\sqrt r}{1+\phi_0(x)/\sqrt s}
\le\sqrt{\ff rs},\quad 0<s\le r,$$
we obtain
\beg{align*}
H_x(s)-H_x(r)&\le2c_2+\ff d2\log\ff rs+2\log\ff{\phi_0(x)+\sqrt r}{\phi_0(x)+\sqrt s}\\
&\le2c_2+\ff{d+2}2\log\ff rs,\ \ 0<s\le r\le 1,\ x\in M^\circ.
\end{align*} This and   \eqref{IU} implies the desired estimate  \eqref{FI}  for some $c\in (0,\infty)$.
 \end{proof}

  We are ready to prove the following result by using  \eqref{WS} and \eqref{FI}.

\beg{lem}\label{FLUX}
There exists a constant $c\in (0,\infty)$  such that  for any $\nu\in\scr P_0$ and any measurable function
$r:M\to [0,1]$:
\beq\label{UE}
\beg{split}
\W_2(\bar\mu_t^\nu,\mu_0)^2&\le\ff ct\int_{M}r(x)^2\nu_0(\d x)+\ff c{t^2}\int_{M}\nu_0(\d x)\int_{r(x)}^{1}s^2q_s^\nu(x)\d s,\qquad t\ge 2s_1.
\end{split}
\end{equation}
The constant $c$ is independent of $\nu$ and $r$.
\end{lem}

\beg{proof} By using $\vv\lor r(x)$ in place of $r(x)$ and letting $\vv\downarrow 0$, we may and do assume that $ \vv_r:=\inf  r>0.$
For any $\vv\in (0, \vv_r)$, let
$$
g_{t,\vv}^\nu:=Q_\vv g_t^\nu=\ff1t\int_\vv^{t+\vv}q_s^\nu\d s.
$$
By  $\mu_0(g_{t}^\nu)=1,$ \eqref{UC} and \eqref{MD2-nu}, we have
$$g_{t,\vv}^\nu\in C_b^1(M),\ \ \ \mu_0(g_{t,\vv}^\nu)=1,\ \ \ g_{t,\vv}^\nu\ge1/4, \ \ \ t\ge 2s_1. $$ So, by \eqref{QE} and the semigroup property,
we obtain
\beg{align*}
\nn(-L_0)^{-1}(g_{t,\vv}^\nu-1)&=\int_0^\infty\nn Q_u(g_{t,\vv}^\nu-1)\d u=\ff 1 t \int_0^\infty \int_{\vv}^{t+\vv}\nn q_{u+v}^\nu\d v \d u\\
&=\ff 1 t \int_\vv^\infty\left(\int_\vv^{(t+\vv)\wedge s}\d v\right) \nn q_s^\nu \d s=\ff1t\int_\vv^\infty\big((s-\vv)\wedge t\big)\nn q_s^\nu\d s\\
&=\ff{U_{r,\vv}^\nu+V_{r,\vv}^\nu}t+\ff1t\int_1^\infty\big((s-\vv)\wedge t\big)\nn q_s^\nu\d s,
\end{align*}
where
$$
\beg{aligned}
U_{r,\vv}^\nu&:=\int_{M}\nu_0(\d x)\int_\vv^{r(x)}(s-\vv)\nn q_s(x,\cdot)\d s,\\
V_{r,\vv}^\nu&:=\int_{M}\nu_0(\d x)\int_{r(x)}^1(s-\vv)\nn q_s(x,\cdot)\d s.
\end{aligned}
$$
By \eqref{WS},    the Cauchy--Schwarz  and Minkowski inequalities, and $g_{t,\vv}^\nu\ge1/4$, we obtain
\beq\label{JT}
\beg{split}
&\W_2(g_{t,\vv}^\nu\mu_0,\mu_0)^2 \le4\mu_0\left(\ff{|\nn(-L_0)^{-1}(g_{t,\vv}^\nu-1)|^2}
{g_{t,\vv}^\nu}\right)\\
&\le\ff{12}{t^2}\mu_0\left(\ff{|U_{r,\vv}^\nu|^2}{g_{t,\vv}^\nu}\right)
+\ff{48}{t^2}\|V_{r,\vv}^\nu\|_{L^2(\mu_0)}^2 +
 \ff{48}{t^2}\left(\int_1^\infty s\|\nn q_s^\nu\|_{L^2(\mu_0)}\d s\right)^2.
\end{split}
\end{equation}
We now estimate the three terms on the right hand side, respectively.

Firstly, by the definition of $g_{t,\vv}^\nu$ and $\vv<r(x)\le1/2<t$,
\beg{align*}
\int_{M}\nu_0(\d x)\int_\vv^{r(x)}q_s(x,y)\d s&\le\int_{M}\nu_0(\d x)\int_\vv^{t+\vv}q_s(x,y)\d s=tg_{t,\vv}^\nu(y).
\end{align*}
This and   the Cauchy--Schwarz inequality yield
\beg{align*}
|U_{r,\vv}^\nu(y)|^2&\le\left(\int_{M}\nu_0(\d x)\int_\vv^{r(x)}q_s(x,y)\d s\right)
 \left(\int_{M}\nu_0(\d x)\int_\vv^{r(x)}(s-\vv)^2\ff{|\nn_yq_s(x,y)|^2}{q_s(x,y)}\d s
\right)\\
&\le tg_{t,\vv}^\nu(y)\int_{M}\nu_0(\d x)\int_\vv^{r(x)}(s-\vv)^2\ff{|\nn_yq_s(x,y)|^2}{q_s(x,y)}\d s,
\end{align*}
so that  \eqref{FI} implies   that for some constant $c_1\in (0,\infty)$
\beq\label{JT1}
\beg{split}
&\mu_0\left(\ff{|U_{r,\vv}^\nu|^2}{g_{t,\vv}^\nu}\right) \le t\int_{M}\nu_0(\d x)\int_\vv^{r(x)}(s-\vv)^2\mu_0\left(\ff{|\nn q_s(x,\cdot)|^2}{q_s(x,\cdot)}\right)\d s\\
&\le t\int_{M}\nu_0(\d x)\int_0^{r(x)}s^2\mu_0\left(\ff{|\nn q_s(x,\cdot)|^2}{q_s(x,\cdot)}\right)\d s
 \le c_1t\int_{M}r(x)^2\nu_0(\d x).
\end{split}
\end{equation}

Next, by the spectral expansion of $q_s$ in \eqref{QE},
$$
\|V_{r,\vv}^\nu\|_{L^2(\mu_0)}^2=\sum_{m=1}^\infty(\ll_m-\ll_0)\left|\int_{r(x)}^1(s-\vv)\e^{-(\ll_m-\ll_0)s}
\int_{M}\psi_m(x)\nu_0(\d x)\d s\right|^2.
$$
Since
$$
\int_0^1(\ll_m-\ll_0)\e^{-(\ll_m-\ll_0)s}\d s\le1,\qquad m\ge1,
$$
 the Cauchy--Schwarz inequality implies
\beq\label{JT2K}
\beg{split}
\|V_{r,\vv}^\nu\|_{L^2(\mu_0)}^2&\le\int_\vv^1(s-\vv)^2\sum_{m=1}^\infty\e^{-(\ll_m-\ll_0)s}
\left|\int_{M}1_{\{r(x)\le s\}}\psi_m(x)\nu_0(\d x)\right|^2\d s\\
&\le\int_\vv^1s^2\int_{M}\int_{M}1_{\{r(x)\le s\}}1_{\{r(z)\le s\}}\big(q_s(x,z)-1\big)\nu_0(\d x)\nu_0(\d z)\d s\\
&\le\int_\vv^1s^2\int_{M}\int_{M}1_{\{r(x)\le s\}}q_s(x,z)\nu_0(\d x)\nu_0(\d z)\d s.
\end{split}
\end{equation}
Combining this with
$$
\int_{M}\int_{M}1_{\{r(x)\le s\}}q_s(x,z)\nu_0(\d x)\nu_0(\d z)=\int_{M}1_{\{r(x)\le s\}}q_s^\nu(x)\nu_0(\d x),
$$ we derive
\beq\label{JT2}
\beg{split}
\|V_{r,\vv}^\nu\|_{L^2(\mu_0)}^2&\le  \int_{M}\nu_0(\d x)\int_{r(x)}^{1}s^2q_s^\nu(x)\d s.
\end{split}
\end{equation}

Moreover,   by \eqref{QE},   \eqref{UC} and  the integration by parts formula, we find a constant $c_2\in (0,\infty)$ such that
\beg{align*}
&\|\nn q_s(x,\cdot)\|_{L^2(\mu_0)}^2 =\sum_{m=1}^\infty(\ll_m-\ll_0)\e^{-2(\ll_m-\ll_0)s}\psi_m(x)^2\\
&\le\sum_{m=1}^\infty\e^{-(\ll_m-\ll_0)(2s-1)}\psi_m(x)^2
 =q_{2s-1}(x,x)-1\le c_2\e^{-2(\ll_1-\ll_0)s},\ \ \  s\ge 1.
\end{align*}
Since $\nu_0$ is a probability measure, by Minkowski's inequality, we obtain
$$
\int_1^\infty s\|\nn q_s^\nu\|_{L^2(\mu_0)}\d s \le\int_1^\infty s\int_{M}\|\nn q_s(x,\cdot)\|_{L^2(\mu_0)}\nu_0(\d x)\d s\\
 \le c_2\int_1^\infty s\e^{-(\ll_1-\ll_0)s}\d s<\infty.$$
 Substituting this, \eqref{JT1} and \eqref{JT2}   into \eqref{JT}, we find a constant $c\in (0,\infty)$ such tat
\beq\label{UPP}
\W_2(g_{t,\vv}^\nu\mu_0,\mu_0)^2\le\ff ct\int_{M}r(x)^2\nu_0(\d x)+\ff c{t^2}\int_{M}\nu_0(\d x)\int_{r(x)}^{1/2}s^2q_s^\nu(x)\d s,\ \ \ \vv\in (0, \vv_r).
\end{equation}

Finally, since $g_{t,\vv}^\nu\mu_0=\bar\mu_t^\nu Q_\vv$, we have
$$
\pi_{t,\vv}^\nu(\d y,\d z):=\bar\mu_t^\nu(\d y)Q_\vv(y,\d z)\in\C\big(\bar\mu_t^\nu,g_{t,\vv}^\nu\mu_0\big),
$$ so that    \eqref{QM} implies
\beg{align*}
\W_2\big(g_{t,\varepsilon}^\nu\mu_0,\bar\mu_t^\nu\big)^2\le\int_{M\times M}
\rr(y,z)^2\pi_{t,\vv}^\nu(\d y,\d z)=\int_{M}Q_\vv\big(\rr(y,\cdot)^2\big)(y)\bar\mu_t^\nu(\d y)\le c_\rr\vv.
\end{align*}
So, by letting $\vv \to 0$ in \eqref{UPP}   we obtain \eqref{UE}.
\end{proof}

\beg{proof}[Proof of the upper bound estimate in Theorem \ref{T1}]
By Lemma \ref{JL},  Lemma \ref{FLUX} and the triangle inequality, there exists a constant $c\in (0,\infty)$  such that
\beq\label{TRI}
\W_2(\mu_t^\nu,\mu_0)^2\le\ff ct\int_Mr(x)^2\nu_0(\d x)+\ff c{t^2}\int_M\nu_0(\d x)
\int_{r(x)}^{1/2}s^2q_s^\nu(x)\d s+\ff c{t^2},\qquad t\ge t_0
\end{equation}
holds for  any  $\nu\in \scr P_0$ and  any measurable $r:M\to [0,1]$.

By \eqref{UC}, we may choose lager $t_0$ such that
$$
\sup_{\nu\in\scr P_0}\sup_{x\in M^\circ}\int_{1/2}^1q_u^\nu(x)\d u<t_0.
$$
For $t\ge t_0$ and $x\in M^\circ$, let
\beq\label{TX-nu}
r(x):=\sup\left\{s\in(0,1):\int_s^1q_u^\nu(x)\d u>t\right\},
\end{equation}
where $\sup\emptyset:=0$.  Then
\beq\label{QEt}
\int_{r(x)}^1q_u^\nu(x)\d u
\le t.  \end{equation}   Combining this with    $$u^2=r(x)^2+2\int_{r(x)}^u s \d s,\quad u>0,$$
we obtain
\beq\label{UQ}
\int_{r(x)}^1u^2q_u^\nu(x)\d u\le tr(x)^2+2\int_{r(x)}^1s\int_s^1q_u^\nu(x)\d u\d s.
\end{equation}
Since
\beq\label{QTU}
\Phi_t(a):=\int_0^a\ff{tu}{t+u}\d u=ta-t^2\log\left(1+\ff at\right),\qquad a\ge0,\  t\ge 1
\end{equation}
satisfies
$$
\Phi'_t(a)=\ff{ta}{t+a}\ge
\begin{cases}
t/2,&a\ge t,\\
a/2,&0\le a\le t,
\end{cases}
$$
we deduce from     \eqref{UQ} and \eqref{QEt} and that
\beg{align*}
tr(x)^2+\int_{r(x)}^{1}s^2q_s^\nu(x)\d s \le2tr(x)^2+2\int_{r(x)}^1s\int_s^1q_u^\nu(x)\d u\d s
 \le8\int_0^1s \Phi'_t\left(\int_s^1q_u^\nu(x)\d u\right)\d s.
\end{align*}
Substituting this into \eqref{TRI}, we find a constant $c'\in (0,\infty)$ such that
\beq\label{W2UN}
\W_2(\mu_t^\nu,\mu_0)^2\le\ff {c'}{t^2}\left[1+\int_0^1s\int_M\Phi_t'\left(\int_s^1q_u^\nu(x)\d u\right)\nu_0(\d x)\d s\right],
\qquad t\ge t_0.
\end{equation}
It remains to   estimate the integral on the right hand side.

Fix $k_1>1$ such that $k_1^2\ge2c_\rr$. By \eqref{QM}, we have
\begin{equation}\label{QBL}
Q_s1_{B(x,k_1\sqrt s)}(x)\ge1-\ff{Q_s(\rr(x,\cdot)^2)(x)}{k_1^2s}
\ge\ff12,\ \ \ s\in (0,1].
\end{equation}
If $\rr(x,y)\le k_1\sqrt s$, then
$$
\rr(z,y)^2\le2\rr(z,x)^2+2k_1^2s,\qquad \phi_0(y)+\sqrt s\le(1+ck_1\|\nn\phi_0\|_\infty)(\phi_0(x)+\sqrt s).
$$
So, by \eqref{QGA},   we   find a constant $k_2\in (0,1)$ such that
\beq\label{QSL}
q_{c_2s}(z,y)\ge k_2 q_s(z,x),\qquad s\in (0, k_2],\  y\in B(x,k_1\sqrt s)\cap M^\circ,\  x,z\in M^\circ.
\end{equation}
By \eqref{UC},
\beq\label{c*}
c_*:= 1+\sup_{\nu\in\scr P_0}\left\|\int_1^\infty(q_u^\nu-1)\d u\right\|_\infty<\infty.
\end{equation}
Then, by the definition of $g^\nu$ in \eqref{V0} and that $q_t$ is the symmetric heat kernel of $Q_t$ with respect to $\mu_0$,  we obtain
\beq\begin{split}\label{QSG}
 Q_sg^\nu+c_*
 &=\int_s^1q_u^\nu\d u+s+c_*-1+\int_1^\infty(q_u^\nu-1)\d u\\
&\ge\int_s^1q_u^\nu\d u+s>0,\qquad s\in (0,1].
\end{split}\end{equation}
This and  \eqref{QSL} imply that for some constant $k\in (1,\infty)$
$$
\int_s^1q_u^\nu(x)\d u \le k(Q_{s/k_2}g^\nu(y)+c_*),\qquad  s\in (0, k_2],\quad y\in B(x,k_1\sqrt s)\cap M^\circ.
$$
Since $\Phi_t'$ is increasing and $\Phi_t'(ka)\le k\Phi_t'(a)$ for $a  \ge0$, this implies
$$
\Phi_t'\left(\int_s^1q_u^\nu(x)\d u\right)\le k\Phi_t'\bigl((Q_{s/k_2}g^\nu(y)+c_*)\bigr),\qquad  s\in (0,k_2],\quad y\in B(x,k_1\sqrt s)\cap M^\circ.
$$
Integrating this inequality over $B(x,k_1\sqrt s)$ against $q_s(x,y)\mu_0(\d y)$ and using \eqref{QBL}, we obtain
\beg{align*}
&\Phi_t'\left(\int_s^1q_u^\nu(x)\d u\right)
 \le 2Q_s1_{B(x,k_1\sqrt s)}(x) \Phi_t'\left(\int_s^1q_u^\nu(x)\d u\right)\\
&\le 2kQ_s\left(1_{B(x,k_1\sqrt s)}\Phi_t'(Q_{s/k_2}g^\nu+c_*)\right)(x)
 \le 2kQ_s\bigl(\Phi_t'(Q_{s/k_2}g^\nu+c_*)\bigr)(x),\qquad s\in (0,k_2].
\end{align*}
Combining this with the following formula implied by \eqref{V0} and the symmetry of $q_s$
$$\nu_0(Q_s f)= \int_{M\times M} q_s(x,y) f(y) \nu_0(\d x) \mu_0(\d y)= \mu_0(q_s^\nu f), \ \ f\ge 0,$$
we arrive at
\beg{equation} \label{Phi'L}
\int_M \Phi_t'\left(\int_s^1q_u^\nu(x)\d u\right)\d\nu_0\le 2 k \mu_0\big(q_s^\nu \Phi_t'(Q_{s/k_2}g^\nu+c_*)\big).\end{equation}

By   \eqref{QGA} and \eqref{QS} we find a constant $k_3\in (1,\infty)$ such that
$$q_s^\nu\le k_3 q_{s/k_2+u}^\nu,\qquad  \quad 0< u\le s\le \ff 1 2 k_2.$$
This together with \eqref{QSG} yields that for some constant  $k_4 \in (1,\infty)$
\begin{equation}\label{sQL}
sq_s^\nu(y)\le   k_3\int_0^1q_{s/k_2+u}^\nu(y)\d u\le k_4(Q_{s/k_2}g^\nu(y)+c_*),\quad  0<s\le \ff 1 2 k_2.
\end{equation}
For $s\in [\ff 1 2 k_2,1]$, \eqref{UC} gives a uniform bound for $\int_s^1q_u^\nu(x)\d u$. Combining this with \eqref{Phi'L},  \eqref{sQL}
and noting that
$a\Phi'_t(a)\le 2\Phi_t(a)$,  we find a constant $k_5\in (1,\infty)$ such that
\beg{equation}\begin{split}\label{Phi'U}
&\int_0^1s\int_M\Phi_t'\left(\int_s^1q_u^\nu(x)\d u\right)
\nu_0(\d x)\d s \le k_5+k_5 \int_0^{\ff 1 2 k_2}s\,\mu_0\bigl(q_s^\nu\Phi_t'(Q_{c_2s}g^\nu+c_*)\bigr)\d s\\
& \le k_5 +2kk_5\int_0^{\ff 1 2 k_2}\mu_0\bigl((Q_{c_2s}g^\nu+c_*)\Phi_t'(Q_{c_2s}g^\nu+c_*)\bigr)\d s \\
& \le k_5+4kk_5\int_0^1\mu_0\bigl(\Phi_t(Q_sg^\nu+c_*)\bigr)\d s.
\end{split}\end{equation}

By $\pp_s Q_s= L_0 Q_s$ and the integration by parts formula, we obtain
\beq\label{DT0}
-\ff{\d}{\d s}\mu_0\bigl(\Phi_t(Q_sg^\nu+c_*)\bigr)=\mu_0\left(\ff{|\nn Q_sg^\nu|^2}{(1+t^{-1}(Q_sg^\nu+c_*))^2}\right),\qquad s>0.
\end{equation}
Since $\mu_0(Q_sg^\nu)=0$, Jensen's inequality for the convex function $\Phi_t$ and $0\le\Phi_t(a)\le ta$ imply
$$
0\le\mu_0\bigl(\Phi_t(Q_sg^\nu+c_*)\bigr)-\Phi_t(c_*)\le tc_*,\qquad s\in (0,1],
$$
while Taylor's expansion  at $c_*$, $\mu_0(Q_sg^\nu)=0$, $\Phi_t''\le1$ and
\eqref{UC} yield that for some constant $k_6\in (1,\infty)$
\beq\label{Phi-D}
0\le\mu_0\bigl(\Phi_t(Q_sg^\nu+c_*)\bigr)-\Phi_t(c_*)\le\ff12\mu_0\bigl((Q_sg^\nu)^2\bigr)
\le k_6\e^{-2(\ll_1-\ll_0)s},\ \ \ s\ge 1.
\end{equation} %Therefore,
Consequently, the function $s\mapsto\mu_0\bigl(\Phi_t(Q_sg^\nu+c_*)\bigr)-\Phi_t(c_*)$ is integrable on $(0,\infty)$, and its product with $s$ tends to zero as   $s\to\infty$.
So,  \eqref{DT0} and the  integration by parts formula imply
\beq\label{DT}
\beg{split}\int_0^\infty s\,\mu_0\left(\ff{|\nn Q_sg^\nu|^2}{(1+t^{-1}(Q_sg^\nu+c_*))^2}\right)\d s=
 \int_0^\infty\left[\mu_0\bigl(\Phi_t(Q_sg^\nu+c_*)\bigr)-\Phi_t(c_*)\right]\d s  <\infty.
\end{split}
\end{equation}
By \eqref{QSG}, we have $Q_sg^\nu\ge-c_*$ for $s>0$. Hence
$$
\bigl|(Q_sg^\nu+c_*)-|Q_sg^\nu|\bigr|\le c_*.
$$
Since $t\ge1$, it follows that
$$
(1+c_*)^{-1}(1+t^{-1}|Q_sg^\nu|)\le1+t^{-1}(Q_sg^\nu+c_*)\le(1+c_*)(1+t^{-1}|Q_sg^\nu|),\qquad s>0.
$$
Combining this with \eqref{DT} and the definition of $\mathcal D_\nu(t)$, we find a constant $k_7\in (1,\infty)$ such that
\beq\label{D-PHI}
\beg{split}
k_7^{-1}\mathcal D_\nu(t)\le\int_0^\infty\left[\mu_0\bigl(\Phi_t(Q_sg^\nu+c_*)\bigr)-\Phi_t(c_*)\right]\d s\le k_7\mathcal D_\nu(t)<\infty,\qquad t\ge1.
\end{split}
\end{equation}
This together with  \eqref{Phi-D} and $\Phi_t(c_*)\le c_*^2/2$  imply
\begin{equation*}
\int_0^1\mu_0\bigl(\Phi_t(Q_sg^\nu+c_*)\bigr)\d s
 \le\int_0^\infty\left[\mu_0\bigl(\Phi_t(Q_sg^\nu+c_*)\bigr)-\Phi_t(c_*)\right]\d s+\ff{c_*^2}2\le k_7  \mathcal D_\nu(t) + c_*^2,
\qquad t\ge1.
\end{equation*}
Combining this with \eqref{W2UN} and \eqref{Phi'U}, we prove the upper bound in \eqref{MAIN-D} for some constant $C\in (1,\infty)$. \end{proof}

\section{Proof of  Theorem \ref{T1}: the lower bound estimate}

For $q_s^\nu$ in \eqref{V0}, let
$$
\rr_t^\nu:=\ff1{t\nu_0(Q_t\phi_0^{-1})}\int_0^t\Big(q_s^\nu Q_{t-s}\phi_0^{-1}-\nu_0(Q_t\phi_0^{-1})
\Big)\d s.
$$
By \cite[(3.4)]{W21}, we have
\beq\label{MX}
\ff{\d\mu_t^\nu}{\d\mu_0}=1+\rr_t^\nu.
\end{equation}
 We first compare $\rr_t^\nu$ with
 \beq\label{TTX}
\tt\rr_t^\nu:=\ff1{t\nu_0(Q_t\phi_0^{-1})}\int_0^\infty\Big(\mu(\phi_0)(q_s^\nu-1)+(Q_s-\mu_0)(\phi_0^{-1})\Big)\d s.
\end{equation}

\beg{lem} Let $s_0$ be in $\eqref{QB}$. There exists a constant $c\in (0,\infty)$  such that
\beq\label{ER}
\mu_0\left(|\rr_t^\nu-\tt\rr_t^\nu|\right)\le c\e^{-t(\ll_1-\ll_0)},\ \ \ t\ge s_0,\ \nu\in \scr P_0.
\end{equation}
\end{lem}

\beg{proof} Let
\beg{align*}
I_s^\nu&:=(q_s^\nu-1)\big(Q_{t-s}\phi_0^{-1}-\mu(\phi_0)\big)+\mu(\phi_0)-\nu_0(Q_t\phi_0^{-1}),\\
J_s^\nu&:=\mu(\phi_0)(q_s^\nu-1)+(Q_s-\mu_0)(\phi_0^{-1}).
\end{align*} We have
$$
\rr_t^\nu-\tt\rr_t^\nu=\ff1{t\nu_0(Q_t\phi_0^{-1})}\bigg(\int_0^t I_s^\nu\d s-\int_t^\infty J_s^\nu\d s\bigg).
$$
By \eqref{QB},
$$
\nu_0(Q_t\phi_0^{-1})\ge\ff12\mu_0(\phi_0^{-1})=\ff12\mu(\phi_0),\ \ \ t\ge s_0.
$$
Therefore,
\beq\label{AS}
\beg{split}
&\mu_0\left(|\rr_t^\nu-\tt\rr_t^\nu|\right)\\
&\le\ff2{t\mu(\phi_0)}\bigg(\int_0^t\|I_s^\nu\|_{L^1(\mu_0)}\d s+\int_t^\infty\|J_s^\nu\|_{L^1(\mu_0)}\d s
\bigg),\ \ \ t\ge s_0.
\end{split}
\end{equation}

Next, by \eqref{UC}, we find  $c_1\in (0,\infty)$ such that
\begin{equation}\begin{split}\label{SQI}
\|q_s^\nu-1\|_\infty&=\sup_{y\in M^\circ}\left|\int_{M}(q_s(x,y)-1)\nu_0(\d x)\right|\\
&\le\sup_{x,y\in M^\circ}|q_s(x,y)-1|\le c_1(1\land s)^{-\ff{d+2}{2}}\e^{-(\ll_1-\ll_0)s},\ \ \ s>0.
\end{split}\end{equation}
Combining this with
 $$\|q_s^\nu-1\|_{L^1(\mu_0)}\le \min\big\{2,\ \|q_s^\nu-1\|_\infty\big\}$$   we find
  $c_2\in (0,\infty)$ such that
\begin{equation}\label{SQI-1}
\|q_s^\nu-1\|_{L^1(\mu_0)}\le c_2\e^{-(\ll_1-\ll_0)s},\qquad s>0.
\end{equation}
By \eqref{QE},
\begin{equation}\label{QL1}
\|Q_s\phi_0^{-1}-\mu(\phi_0)\|_{L^1(\mu_0)}\le\|Q_s\phi_0^{-1}-\mu(\phi_0)\|_{L^2(\mu_0)}\le\e^{-(\ll_1-\ll_0)s},\ \ s>0.
\end{equation}
By \eqref{UC}, $\mu_0(q_s(x,\cdot))=1$ and $\mu_0(q_s(x,\cdot)^2)=q_{2s}(x,x)$,   we find  $c_3\in (0,\infty)$ such that
\beg{align*}
\int_M|q_s(x,y)-1|^2\mu_0(\d y)=q_{2s}(x,x)-1\le c_3(1\land s)^{-\ff{d+2}{2}}\e^{-2(\ll_1-\ll_0)s},\qquad x\in M^\circ,\quad s>0.
\end{align*}
Combining this with the Cauchy--Schwarz inequality  and $\mu_0(\phi_0^{-2})=1$,   we obtain
\begin{equation}\label{BQI}\beg{split}
&\big|\mu(\phi_0)-\nu_0(Q_s\phi_0^{-1})\big|\le\|Q_s\phi_0^{-1}-\mu(\phi_0)\|_\infty\\
&= \sup_{x\in M^\circ} \bigg|\int_M \big(q_s(x,y)-1\big)\phi_0^{-1}(y) \mu_0(\d y)\bigg|\le
c_3^{1/2}(1\land s)^{-\ff{d+2}{4}}\e^{-(\ll_1-\ll_0)s},\qquad s>0.
\end{split}\end{equation}
By \eqref{SQI}-\eqref{BQI}, we find   $c_4\in (0,\infty)$  such that
\beq\label{I1-}
\beg{split}
&\int_0^t\|I_s^\nu\|_{L^1(\mu_0)}\d s
 \le\int_0^{\ff t2}\|q_s^\nu-1\|_{L^1(\mu_0)}\|Q_{t-s}\phi_0^{-1}-\mu(\phi_0)\|_{L^\infty(\mu_0)}\d s\\
&\quad+\int_{\ff t2}^t\|q_s^\nu-1\|_{L^\infty(\mu_0)}\|Q_{t-s}\phi_0^{-1}-\mu(\phi_0)\|_{L^1(\mu_0)}
\d s +t\big|\mu(\phi_0)-\nu_0(Q_t\phi_0^{-1})\big|\\
&\le c_4\int_0^t\e^{-s(\ll_1-\ll_0)}\e^{-(t-s)(\ll_1-\ll_0)}\d s+c_4t\e^{-t(\ll_1-\ll_0)}
 \le 2c_4t\e^{-t(\ll_1-\ll_0)},\ \ \ t\ge s_0,\ \nu\in \scr P_0.
\end{split}
\end{equation}

Finally, by using \eqref{SQI-1} and \eqref{QL1} again, we find  $c_5\in (0,\infty)$ such that
\beg{align*}
\int_t^\infty\|J_s^\nu\|_{L^1(\mu_0)}\d s&\le\int_t^\infty\Big(\mu(\phi_0)\|q_s^\nu-1\|_{L^1(\mu_0)}+\|(Q_s-\mu_0)(\phi_0^{-1})\|_{L^1(\mu_0)}
\Big)\d s\\
&\le c_5\e^{-t(\ll_1-\ll_0)},\ \ \ t\ge s_0,\ \nu\in \scr P_0.
\end{align*}
Combining this with \eqref{AS} and \eqref{I1-}, we finish the proof.
\end{proof}
\black

\beg{lem}\label{KG-L}
There exist   $c,k\in (1,\infty)$  such that,  for any $s\in (0,1]$ and bounded nonnegative measurable function $f$ on $M$,
\beq\label{QDG-L}
\|\nn Q_sf\|_\infty\le c  s^{-1/2}\|f\|_\infty,
\end{equation}
\beq\label{QHD-L}
|\nn Q_sf|(z)\le cs^{-1/2}Q_{k s}f(y),\qquad  y,z\in M^\circ,\ \rr(y,z)\le2\sqrt s.
\end{equation}
\end{lem}

\beg{proof} We will apply the following Bismut formula
\cite[Lemma 2.2]{WJEMS} with $\gamma(u)=u/s$:
\beq\label{BSM}
\nn Q_sf(x)=\ff1{\sqrt2}\E^x_Q\left[f(Y_s)\int_0^s\ff1s R_u^*\d B_u\right],\qquad x\in M^\circ,\end{equation}
where $\E_Q^x$ is the expectation taken for the diffusion process $(Y_s)_{s\ge 0}$ generated by $L_0$ with initial value $x$, $B_s$ is the underlying $d$-dimensional Brownian motion,  $R_s^*$ is the transposition of the curvature flow $(R_s)_{s\ge 0}$ on $\R^d\otimes \R^d$ solving the equation
\beq\label{RO-L}
\pp_s R_s=-\operatorname{Ric}_{L_0}^{\#}(U_s)R_s,\qquad R_0=I,\ \ s\ge 0,
\end{equation}
for $(U_s)_{s\ge 0}$   the horizontal lift   of   $(Y_s)_{s\ge 0}$ to the frame bundle $O(M^\circ)$,  and $\operatorname{Ric}_{L_0}^{\#}(U_s)\in \R^d\otimes\R^d$ is defined via
$$
\big\langle\operatorname{Ric}_{L_0}^{\#}(U_s)\xi,\eta\big\rangle_{\R^d}:=
({\rm Ric}-{\rm Hess}_{V+2\log \phi_0})(U_s\xi, U_s\eta),\ \ \ \xi,\eta\in\R^d.$$
  Since $V,\phi_0\in C^2(M)$ and $\phi_0>0$ in $M^\circ$, there exists a constant $c_0\in (0,\infty)$ such that
\beq\label{RC-L}
\operatorname{Ric}_{L_0}
 =\operatorname{Ric}-\operatorname{Hess}V-2\phi_0^{-1}\operatorname{Hess}\phi_0+2\phi_0^{-2}\d\phi_0\otimes\d\phi_0
\ge-c_0\phi_0^{-1}.
\end{equation}

(a) Proof of \eqref{QDG-L}.
 By \eqref{DU},
$$
P_s^D1(y)\le c\left(1\land\ff{\rr_\pp(y)}{\sqrt s}\right),\ \ \ s\in (0,1].
$$
By \eqref{QK}, \eqref{DU}   and \eqref{KK}, we find $c_1\in (0,\infty)$ such that
\beq\label{BI-L}
Q_s\phi_0^{-1}(y)=\ff{\e^{\ll_0s}P_s^D1(y)}{\phi_0(y)}\le c_1s^{-1/2},\qquad y\in M^\circ,\quad 0<s\le1.
\end{equation}
Taking $a:= (8c_0c_1)^{-2},$ and applying the Markov property of $Y_s$, we obtain
\beg{align*}
&\E_Q^x\left[2c_0\int_u^a\phi_0^{-1}(Y_v)\d v\,\middle|\,Y_w,\ 0\le w\le u\right]
 =2c_0\int_0^{a-u}Q_v(\phi_0^{-1})(Y_u)\d v\\
 &\le   4c_0c_1  \sqrt{a-u}\le  \ff12,\ \ \ x\in M^\circ,\ u\in [0,a].
\end{align*}
So, for any $n\ge2$,
\beg{align*}
&\sup_{x\in M^\circ}\E_Q^x\left[\left(2c_0\int_0^a\phi_0^{-1}(Y_u)\d u\right)^n\right]
 =n!\sup_{x\in M^\circ}\int_{0<u_1<\cdots<u_n<a}
\E_Q^x\left[\prod_{j=1}^n2c_0\phi_0^{-1}(Y_{u_j})\right]\d u_1\cdots\d u_n\\
&\quad\le\ff{n!}2\sup_{x\in M^\circ}\int_{0<u_1<\cdots<u_{n-1}<a}\E_Q^x\left[\prod_{j=1}^{n-1}2c_0\phi_0^{-1}(Y_{u_j})
\right]\d u_1\cdots\d u_{n-1}\\
&\quad=\ff n2\sup_{x\in M^\circ}\E_Q^x\left[\left(2c_0\int_0^a\phi_0^{-1}(Y_u)\d u\right)^{n-1}\right]\le\cdots\le  n! 2^{-n}.
\end{align*}
Consequently,
$$
\sup_{x\in M^\circ}\E_Q^x\left[\exp\left(2c_0\int_0^a\phi_0^{-1}(Y_u)\d u\right)\right]\le\sum_{n=0}^\infty2^{-n}=2.
$$
Letting $n_0\in \mathbb N$ such that $n_0a\ge 1$, this together with the Markov property of $Y_s$ yields
\beq\label{KE-L}
 \sup_{x\in M^\circ}\E_Q^x\left[\exp\left(2c_0\int_0^1\phi_0^{-1}(Y_u)\d u\right)\right]
\le2^{n_0}<\infty.
\end{equation}
By \eqref{RO-L},
$$
\beg{aligned}
\ff{\d}{\d t}|R_t\xi|^2&=-2\left\langle\operatorname{Ric}_{L_0}^{\#}(U_t)R_t\xi,R_t\xi
\right\rangle
 =-2\operatorname{Ric}_{L_0}(Y_t)\big(U_tR_t\xi,U_tR_t\xi\big)\\
&\le2c_0\phi_0^{-1}(Y_t)|R_t\xi|^2,\ \ \ t\ge 0,\ \xi\in\R^d,
\end{aligned}
$$
where the last inequality follows from \eqref{RC-L} and $|U_tR_t\xi|=|R_t\xi|$, since $U_t$ is an orthonormal frame. Then, by the Gronwall's inequality,
$$
\|R_t\|^2\le\exp\left(2c_0\int_0^t\phi_0^{-1}(Y_u)\d u\right),
$$
where $\|R_t\|:=\sup_{|v|=1}|R_tv|$ denotes the operator norm. This together with  \eqref{KE-L} implies
$$
 \sup_{x\in M^\circ}\E_Q^x\bigg[\sup_{0<s\le1}\|R_s\|^2\bigg]<\infty,
$$
Combining this with \eqref{BSM},
  the Cauchy-Schwarz inequality and It\^o's isometry,  we find a constant $c\in (0,\infty)$ such that
\beq\begin{split}\label{QGP-L}
&|\nn Q_sf(x)|^2 \le\ff12Q_s(f^2)(x)\E_Q^x\left|\int_0^s\ff1sR_u^*\d B_u\right|^2\\
&\le\ff d{2s^2}Q_s(f^2)(x)\int_0^s\E_Q^x\|R_u\|^2\d u\le\ff csQ_s(f^2)(x),\ \ x\in M^\circ,\ s\in (0,1].
\end{split}\end{equation}
Since $Q_s(f^2)\le\|f\|_\infty^2$, this proves \eqref{QDG-L}.

(b) Proof of  \eqref{QHD-L}. Let $c_1\in (1,\infty)$ be in \eqref{QGA}, and choose $b:= 8 c_1^2>1.$ Then \eqref{QGA} implies
$$
\ff{q_s(x,z)}{q_{bs}(x,y)}\le c_1^2b^{d/2}\ff{(\rr_\pp(x)+\sqrt{bs})(\rr_\pp(y)+\sqrt{bs})}
{(\rr_\pp(x)+\sqrt s)(\rr_\pp(z)+\sqrt s)}\exp\left(-\ff{\rr(x,z)^2}{c_1s}+\ff{c_1\rr(x,y)^2}{bs}\right),\ \ s\in (0,b^{-1}].
$$
Since $b\ge 1$ and the triangle inequality
$ |\rr_\pp(y)-\rr_\pp(z)|\le\rr(y,z), $
then for any $s>0$,
\beg{align*}
\ff{(\rr_\pp(x)+\sqrt{bs})(\rr_\pp(y)+\sqrt{bs})}{(\rr_\pp(x)+\sqrt s)(\rr_\pp(z)+\sqrt s)}\le
\ff{(\rr_\pp(x)+\sqrt{bs})(\rr_\pp(z)+\rr(y,z)+\sqrt{bs})}{(\rr_\pp(x)+\sqrt s)(\rr_\pp(z)+\sqrt s)}
\le b\left(1+\ff{\rr(y,z)}{\sqrt s}\right).
\end{align*}
Moreover, by the triangle inequality and $b\ge2c_1^2$, we have
\beg{align*}
-\ff{\rr(x,z)^2}{c_1s}+\ff{c_1\rr(x,y)^2}{bs}\le-\left(\ff1{c_1}-\ff{2c_1}b\right)
\ff{\rr(x,z)^2}s+\ff{2c_1\rr(y,z)^2}{bs}\le\ff{2c_1\rr(y,z)^2}{bs}.
\end{align*}
Combining these estimates, we find $c_2\in (0,\infty)$ such that
\beq\label{KR-U}
q_s(x,z)\le c_2 q_{bs}(x,y)\left(1+\ff{\rr(y,z)}{\sqrt s}\right)\exp\left(\ff{2c_1\rr(y,z)^2}{bs}\right),\qquad 0<s\le b^{-1},
\end{equation}
so that
\beq\begin{split}\label{QSM}
&Q_s\big(q_s(x,\cdot)^2\big)(y) =\int_Mq_s(y,z)q_s(x,z)^2\mu_0(\d z)\\
&\le c_2 q_{bs}(x,y)^2\int_Mq_s(y,z)\left(1+\ff{\rr(y,z)}{\sqrt s}\right)^2\exp\left(\ff{4c_1\rr(y,z)^2}{bs}\right)\mu_0(\d z),\ \ 0<s\le b^{-1}.
\end{split}\end{equation}
Noting that \eqref{KK} and the triangle inequality  imply
$$\phi_0(z)^2 \le \kk^2 \rr_\pp(z)^2\le \kk^2\rr_\pp(z)(\rr_\pp(y)+ \rr(y,z)),$$
by \eqref{QGA}, \eqref{KK}, the boundedness of $V$,
$\mu_0(\d z)=\phi_0(z)^2\e^{V(z)}\d z$ and $b= 8c_1^2$, we find $c_3\in (0,\infty)$ such that
and the triangle inequality, we obtain, for $y,z\in M^\circ$ and $0<s\le b^{-1}$,
\beg{align*}
&q_s(y,z)\phi_0(z)^2\e^{V(z)}\exp\left(\ff{4c_1\rr(y,z)^2}{bs}\right)\\
&\quad\le c_3s^{-d/2}\ff{\rr_\pp(z)(\rr_\pp(y)+\rr(y,z))}{(\rr_\pp(y)+\sqrt s)(\rr_\pp(z)+\sqrt s)}
\exp\left[-\left(\ff1{c_1}-\ff{4c_1}b\right)\ff{\rr(y,z)^2}s
\right]\\
&\quad\le c_3s^{-d/2}\left(1+\ff{\rr(y,z)}{\sqrt s}\right)\exp\left(-\ff{\rr(y,z)^2}{2c_1s}\right).
\end{align*}
Combining this with $\operatorname{vol}(B(y,r))\le cr^d$ for some $c\in (0,\infty)$ and all $r>0$, we obtain
\beg{align*}
&\int_Mq_s(y,z)
\left(1+\ff{\rr(y,z)}{\sqrt s}\right)^2\exp\left(\ff{4c_1\rr(y,z)^2}{bs}\right)\mu_0(\d z)\\
&\quad\le c_3s^{-d/2}\sum_{j=0}^\infty \int_{\{j\ss s\le \rr(y,z)<(j+1)\ss s\}}\left(1+\ff{\rr(y,z)}{\sqrt s}\right)^3
\exp\left(-\ff{\rr(y,z)^2}{2c_1s}\right)\d z\\
&\quad\le c_3s^{-d/2}\sum_{j=0}^\infty(j+2)^3\e^{-j^2/(2c_1)}\operatorname{vol}\big(B(y,(j+1)\sqrt s)\big)\\
&\quad\le c_3c\sum_{j=0}^\infty(j+2)^{d+3}\e^{-j^2/(2c_1)}<\infty,\ \ y\in M^\circ,\ s\in (0,b^{-1}].
\end{align*}
By combining this with \eqref{QSM}, we find $c_4\in (0,\infty)$ such that
\beq\label{K2-L}
\left[Q_s\big(q_s(x,\cdot)^2\big)(y)\right]^{1/2}\le c_4 q_{bs}(x,y),\qquad x,y\in M^\circ,\quad 0<s\le b^{-1}.
\end{equation}

For bounded $f\ge0$, by \eqref{QGP-L}, Minkowski's inequality, and \eqref{K2-L}, we find
\beg{align*}
|\nn Q_{2s}f|(y)&\le cs^{-1/2}\left[Q_s\big((Q_sf)^2\big)(y)\right]^{1/2}
 \le c_2s^{-1/2}\int_M f(x)\left[Q_s\big(q_s(x,\cdot)^2\big)(y)\right]^{1/2}
\mu_0(\d x)\\
&\le c_2c_4s^{-1/2}Q_{bs}f(y),\qquad 0<s\le b^{-1},\ y\in M^\circ.
\end{align*}
Thus, for some $c_5\in (0,\infty)$ and any bounded $f\ge 0$,
\beq\label{QDP-L}
|\nn Q_sf| \le c_5s^{-1/2}Q_{bs/2}f,\qquad 0<s\le2b^{-1}.
\end{equation}

Let $k:=b^2/2$. By \eqref{KR-U} with $s$ replaced by $bs/2$, we find $c_6\in (0,\infty)$ such that
\beg{align*}
&q_{bs/2}(x,z)\le c_2 q_{\kk s}(x,y)\left(1+\sqrt{\ff2b}\,\ff{\rr(y,z)}{\sqrt s}\right)
\exp\left(\ff{4c_1\rr(y,z)^2}{b^2s}\right)\\
&\le c_6 q_{\kk s}(x,y),\ \ \ s\in (0,k^{-1}],\ y\in B\big(z,2\ss s\big).
\end{align*}
This and  \eqref{QDP-L} imply
\beq\label{9*}
|\nn Q_sf|(z)\le c_5s^{-1/2}Q_{bs/2}f(z)\le c_5c_6s^{-1/2}Q_{k s}f(y),\ \  s\in (0,k^{-1}],\ y\in B\big(z,2\ss s\big).\end{equation}

Now let $k^{-1}\le s\le1$. By \eqref{QGA} and the compactness of $M$, we find
$c_7\in (1,\infty)$ such that
\beq\label{HKUL}
c_7^{-1}\le q_s(y,z)\le c_7,\qquad y,z\in M^\circ,\quad s\ge(2k)^{-1}.
\end{equation}
Applying \eqref{QDG-L} at time $(2k)^{-1}$ and then using \eqref{HKUL}, we find $c_8,c_9\in (0,\infty)$ such that
\beg{align*}
&|\nn Q_sf|(z) \le\left\|\nn Q_{(2k)^{-1}}\bigl(Q_{s-(2k)^{-1}}f\bigr)\right\|_\infty
 \le c_8\left\|Q_{s-(2k)^{-1}}f\right\|_\infty \\
 &\le c_7c_8\mu_0(f)\le c_9s^{-1/2}Q_{k s}f(y),\qquad y,z\in M^\circ,\ s\in [k^{-1},1].
\end{align*} This and \eqref{9*} imply \eqref{QHD-L} for some constant $c\in (0,\infty)$.
\end{proof}

\beg{proof}[Proof of the lower bound estimate in Theorem \ref{T1}] Let $t_0$ be in Lemma \ref{JL}, $g^\nu$ in \eqref{G-REP} and  $\Phi_t$  in \eqref{QTU}. We have
 $\Phi_t'(a)= \ff{at}{a+t}\in [0,t\land a]$ for $a\ge 0$,  and $\mu_0(g^\nu)=0$, so that
 \beq\label{SM}\beg{split}
&f_t^\nu:=\int_0^1Q_s\bigl(\Phi_t'(g^\nu+c_*)\bigr)\d s\in [0,t],\\
&\mu_0(f_t^\nu)=\mu_0\bigl(\Phi_t'(g^\nu+c_*)\bigr)\le\mu_0(g^\nu+c_*)=c_*.\end{split}
\end{equation} We claim that for some constant $K\in (1,\infty)$
\beq\label{FD-L}
\beg{split}
 |f_t^\nu(y)-f_t^\nu(z)|^2&\le K \rr(y,z)^2\bigl[\Phi_t'(g^\nu(y)+c_*)f_t^\nu(y)+\Phi_t'(g^\nu(z)+c_*)f_t^\nu(z)\bigr],\\
 &\qquad \ \ y,z\in M^\circ,\ t\ge t_0,\ \nu\in \scr P_0,
\end{split}
\end{equation}
which, together with the optimal coupling $\pi\in \C(\mu_t^\nu,\mu_0), \d\mu_t^\nu= (1+\rr_t^\nu)\d\mu_0$ and the Cauchy-Schwarz inequality, yields
\beq\label{W2L-}\beg{split}
&\left|\mu_0(f_t^\nu\rr_t^\nu)\right|
 =\left|\int_{M\times M}\bigl(f_t^\nu(y)-f_t^\nu(z)\bigr)\pi(\d y,\d z)\right|\\
&\le K\left(\int_{M\times M}\rr(y,z)^2\pi(\d y,\d z)\right)^{1/2}
 \left[\mu_t^\nu\bigl(\Phi_t'(g^\nu+c_*)f_t^\nu\bigr)
+\mu_0\bigl(\Phi_t'(g^\nu+c_*)f_t^\nu\bigr)\right]^{1/2}\\
&\le K\W_2(\mu_t^\nu,\mu_0)\sqrt{\mu_0\bigl((2+\rho_t^\nu)\Phi'_t(g^\nu+c_*)f_t^\nu\bigr)},\ \ t\ge t_0,\ \nu\in \scr P_0.
\end{split}\end{equation}
In the following, we first prove \eqref{FD-L}, then estimate  $\left|\mu_0(f_t^\nu\rr_t^\nu)\right|$ and $\mu_0\bigl((2+\rho_t^\nu)\Phi'_t(g^\nu+c_*)f_t^\nu\bigr)$ respectively.

(a) Proof of \eqref{FD-L}.   Let $k>1$ be the constant in \eqref{QHD-L}.  By \eqref{QSG},
$$
g^\nu+c_*\ge c^{-1},\quad \Phi_t'(g^\nu+c_*)\ge\ff{g^\nu+c_*}{1+g^\nu+c_*}\ge c^{-1},\ \ \ t\ge t_0.
$$
From \eqref{QTU}, we see that   $\Phi_t'$ is increasing, concave and $1$-Lipschitz. So, by Jensen's inequality, \eqref{G-REP} and the semigroup property of $Q_t$,
we find a constant $d_1\in (0,\infty)$ such that
\beq\label{QF-L}
\beg{split}
&Q_u\bigl(\Phi_t'(g^\nu+c_*)\bigr) \le\Phi_t'(Q_ug^\nu+c_*)\\
 &=\Phi_t'\left(g^\nu+c_*+u-\int_0^uq_v^\nu\d v\right)
 \le\Phi_t'(g^\nu+c_*+u) \\
&\le\Phi_t'(g^\nu+c_*)+u\le d_1\Phi_t'(g^\nu+c_*),\qquad u\in (0,k],\ \nu\in \scr P_0.
\end{split}
\end{equation}
For fixed $y,z\in M^\circ$, let $r:=\rr(y,z)$.
  By the Cauchy-Schwarz inequality, \eqref{QF-L}
and \eqref{SM}, we obtain
\beq\label{ST}
\beg{split}
&\left[\int_0^{r^2\wedge1}\left|Q_s\bigl(\Phi_t'(g^\nu+c_*)\bigr)(y)-Q_s\bigl(\Phi_t'(g^\nu+c_*)\bigr)(z)\right|\d s\right]^2\\
& \le (1\land r^2) \int_0^{1\land r^2}2\left[Q_s\bigl(\Phi_t'(g^\nu+c_*)\bigr)(y)^2+Q_s\bigl(\Phi_t'(g^\nu+c_*)\bigr)(z)^2
\right]\d s\\
& \le 2d_1r^2\left[\Phi_t'(g^\nu(y)+c_*)f_t^\nu(y)+\Phi_t'(g^\nu(z)+c_*)f_t^\nu(z)\right].
\end{split}
\end{equation}
Thus, if $r:=\rr(y,z)\ge 1$,   this and \eqref{SM} imply  \eqref{FD-L} for some constant $K$.

Next, assume $r:=\rr(y,z)<1$. Choosing a  smooth curve $\gamma: [0,1]\to M^\circ$ such that
$$
\gamma(0)=y,\qquad\gamma(1)=z,\qquad\int_0^1|\dot\gamma(u)|\d u\le2r.
$$
Then
$$
\rr(y,\gamma(u))\le\int_0^u|\dot\gamma(v)|\d v\le2r\le2\sqrt s,\ \  \ u\in[0,1],\ s\in [r^2,1].
$$
Applying \eqref{QHD-L} along $\gamma$, we find $k,d_2\in (1,\infty)$ such that
\beg{align*}
&\left|Q_s\bigl(\Phi_t'(g^\nu+c_*)\bigr)(y)-Q_s\bigl(\Phi_t'(g^\nu+c_*)\bigr)(z)
\right| \\
& \le\int_0^1\left|\nn Q_s\bigl(\Phi_t'(g^\nu+c_*)\bigr)(\gamma(u))\right| |\dot\gamma(u)|\d u\\
& \le d_2rs^{-1/2}Q_{k s}\bigl(\Phi_t'(g^\nu+c_*)\bigr)(y),\qquad s\in [r^2,1].
\end{align*}
Then
\beq\begin{split}\label{Diff-S}
&\int_{r^2}^1\left|Q_s\bigl(\Phi_t'(g^\nu+c_*)\bigr)(y)-Q_s\bigl(\Phi_t'(g^\nu+c_*)\bigr)(z)
\right|\d s\\
&\quad\le d_2 k^{-1/2}r\int_{k r^2}^k u^{-1/2}Q_u\bigl(\Phi_t'(g^\nu+c_*)\bigr)(y)\d u.
\end{split}\end{equation}
By \eqref{HKUL} and \eqref{SM},
$$
\beg{gathered}\int_0^k Q_s\bigl(\Phi_t'(g^\nu+c_*)\bigr)\d s\le f_t^\nu+c_7\mu_0\bigl(\Phi_t'(g^\nu+c_*)\bigr),\\
 c_7^{-1}\mu_0\bigl(\Phi_t'(g^\nu+c_*)\bigr)\le \int_{1/2}^1Q_s\bigl(\Phi_t'(g^\nu+c_*)\bigr)\d s\le f_t^\nu.
\end{gathered}
$$
Combining this with \eqref{QF-L}, we find $d_3\in (1,\infty)$ such that
\beg{align*}
&\left[\int_0^k s^{-1/2}Q_s\bigl(\Phi_t'(g^\nu+c_*)\bigr)\d s
\right]^2\\
&\quad=2\int_0^k s^{-1/2}Q_s\bigl(\Phi_t'(g^\nu+c_*)\bigr)\d s\int_0^s u^{-1/2}Q_u\bigl(\Phi_t'(g^\nu+c_*)\bigr)\d u \\
&\quad\le 4d_1\Phi_t'(g^\nu+c_*)\int_0^k Q_s\bigl(\Phi_t'(g^\nu+c_*)\bigr)\d s
 \le d_3\Phi_t'(g^\nu+c_*)f_t^\nu.
\end{align*}
By the definition of $f_t^\nu$ in \eqref{SM}, this together with  \eqref{ST} and \eqref{Diff-S} implies  \eqref{FD-L}.

(b) Estimate   $\mu_0\bigl((2+\rho_t^\nu)\Phi'_t(g^\nu+c_*)f_t^\nu\bigr)$.
We claim that there exists $C\in (1,\infty)$ such that
\beq\label{DTT} C^{-1} \mathcal D_\nu(t)-C\le  \mu_0\big((g^\nu+c_*)f_t^\nu\big)\le C+C\mathcal D_\nu(t),\ \ \ t\ge t_0,\ \nu\in\scr P_0.\end{equation}
By  \eqref{Phi-D}, \eqref{D-PHI}  and $\Phi_t(a)\le a\Phi_t'(a)$ for $a\ge 0$, we find $c_0\in (0,\infty)$ such that
\beq\label{ESS}  \begin{split}
\mathcal D_\nu(t)&\le c_0+c_0\int_0^1\mu_0\bigl(\Phi_t(Q_sg^\nu+c_*)\bigr)\d s\\
&\le c_0+c_0\int_0^1\mu_0\bigl((Q_sg^\nu+c_*)\Phi_t'(Q_s g^\nu+c_*)\bigr)\d s, \qquad t\ge t_0,\ \nu\in \scr P_0.
\end{split}\end{equation}
By \eqref{V0} and the semigroup property of $Q_s$, we obtain
\begin{equation*}
Q_sg^\nu=\int_s^\infty(q_u^\nu-1)\d u=g^\nu+s-\int_0^sq_u^\nu\d u\le g^\nu+s,\qquad 0<s\le1.
\end{equation*}
Since  $ \Phi_t''(a)=\ff{t^2}{(t+a)^2}\in (0,1]$ for $a\ge 0$,    we have
$$
\Phi_t'(Q_sg^\nu+c_*)\le\Phi_t'(g^\nu+c_*+s)\le\Phi_t'(g^\nu+c_*)+s.
$$
Combining this with \eqref{ESS} and    the symmetry of $Q_s$ in $L^2(\mu_0)$, we find $c_1\in (1,\infty)$ such that
  \beq\begin{split}\label{D-L1}
\mathcal D_\nu(t)
&\le c_0+c_0\int_0^1\mu_0\bigl((g^\nu+c_*) Q_s\Phi_t'(Q_sg^\nu+c_*)\bigr)\d s\\
&=c_1+c_1\mu_0\bigl((g^\nu+c_*)f_t^\nu\bigr),\qquad t\ge t_0,\ \nu\in \scr P_0.
\end{split}\end{equation}
So, the first inequality in \eqref{DTT} holds for some constant $C\in (1,\infty).$

On the other hand,  by \eqref{QTU},  we have $\Phi_t'''(a)\le 0, \Phi_t(a)\le \ff 1 2 a^2$ and
%$$a\Phi_t''(a)= \ff{at^2}{a+t)^2}\le \ff{at}{a+t}=\Ph_t'(a),\ \ \ a\ge 0,$$ so that
$a\Phi_t'(a)\le 2 \Phi_t(a)$ for $a\ge 0$. Combining this with the symmetric Markov semigroup property of $Q_s$ in $L^2(\mu_0)$,  Jensen's inequality, and \eqref{D-PHI},
we obtain
\beq \label{AT-L2}\begin{split}
&\mu_0\bigl((g^\nu+c_*)f_t^\nu\bigr) =2\int_0^{1/2}\mu_0\left((Q_sg^\nu+c_*)Q_s\bigl(\Phi_t'(g^\nu+c_*)\bigr)\right)\d s\\
&\le2\int_0^{1/2}\mu_0\left((Q_sg^\nu+c_*)\Phi_t'(Q_sg^\nu+c_*)\right)\d s
 \le4\int_0^{1/2}\mu_0\bigl(\Phi_t(Q_sg^\nu+c_*)\bigr)\d s\\
&\le 4C\mathcal D_\nu(t)+2\Phi_t(c_*)
 \le 4C \mathcal D_\nu(t) + c_*^2,\ \ \ t\ge t_0,\ \nu\in \scr P_0.
\end{split}\end{equation} So, the second inequality  in \eqref{DTT} holds for some constant $C\in (1,\infty).$

Finally, by \eqref{UC} and \eqref{BI-L}, we find $c_1\in (1,\infty)$ such that
$$
\|Q_u\phi_0^{-1}\|_\infty\le c_1(1+u^{-1/2})\ \text{for}\  u>0,\qquad\|q_s^\nu\|_\infty\le c_1 \ \text{for}\  s\ge1.
$$ Combining this with   \eqref{MX}, \eqref{QB} and \eqref{QSG}, we find $c_2, c_3,c_4\in (1,\infty)$ such that
\beg{align*}
1+\rr_t^\nu&=\ff1{t\nu_0(Q_t\phi_0^{-1})}\int_0^tq_s^\nu Q_{t-s}\phi_0^{-1}\d s
 \le\ff {c_2}t\int_0^{t/2}q_s^\nu\d s+\ff {c_2}t\int_{t/2}^tQ_{t-s}\phi_0^{-1}\d s\\
&\le\ff {c_3}t(g^\nu+c_*+t)+\ff {c_3}t\int_0^{t/2}(1+u^{-1/2})\d u \\
& \le c_4\left(1+\ff{g^\nu+c_*}t\right),\ \ \ t\ge t_0,\ \nu\in \scr P_0.
\end{align*}
Noting that $(1+a)\Phi'_t(a)=a$ for $a>0$, this  and \eqref{DTT} imply that  for some constant $C'\in (1,\infty) $
 \beq\label{ES*} \mu_0\bigl((2+\rho_t^\nu)\Phi'_t(g^\nu+c_*)f_t^\nu\bigr)\le  C'+C'\mathcal D_\nu(t),\ \ \ t\ge t_0,\ \nu\in \scr P_0.\end{equation}

(c)  Estimate $\left|\mu_0(f_t^\nu\rr_t^\nu)\right|$.  By \eqref{TTX} and \eqref{V0}, we have
$$
t\nu_0(Q_t\phi_0^{-1})\tilde\rho_t^\nu-\mu(\phi_0)(g^\nu+c_*)=\int_0^\infty\bigl(Q_s\phi_0^{-1}-\mu(\phi_0)\bigr)\d s-c_*\mu(\phi_0).
$$
Combining this with  \eqref{BI-L} for $0<s\le1$ and \eqref{BQI} for $s\ge1$, we find $k_1\in (0,\infty)$ such that
\beg{align*}
&\left\|t\nu_0(Q_t\phi_0^{-1})\tilde\rho_t^\nu-\mu(\phi_0)(g^\nu+c_*)\right\|_\infty\\
&\quad\le k_1+k_1  \int_0^1(1+s^{-1/2})\d s+k_1\int_1^\infty\e^{-(\ll_1-\ll_0)s}\d s=:k_2<\infty.
\end{align*}
So, by \eqref{QB} and \eqref{SM},  we find $k_3\in (1,\infty)$ such that
\beg{align*}
t\mu_0(f_t^\nu\tilde\rho_t^\nu)&\ge\ff{\mu(\phi_0)}{\nu_0(Q_t\phi_0^{-1})}\mu_0\bigl((g^\nu+c_*)f_t^\nu\bigr)-\ff{k_2\mu_0(f_t^\nu)}{\nu_0(Q_t\phi_0^{-1})}\\
&\ge k_3^{-1}\mu_0\bigl((g^\nu+c_*)f_t^\nu\bigr)-k_3,\ \ \ t\ge t_0,\ \nu\in \scr P_0.
\end{align*}
 This together with \eqref{D-L1} and \eqref{ER} yields that for some $k_4\in (1,\infty)$
$$
\mathcal D_\nu(t)\le k_4+k_4t|\mu_0(f_t^\nu\rho_t^\nu)|,\qquad t\ge t_0,\ \nu\in \scr P_0.$$
Combining this with \eqref{W2L-} and \eqref{ES*},    we find $k_5\in (1,\infty)$ such that
\beg{align*}
&\mathcal D_\nu(t) \le k_4+k_4Kt \W_2(\mu_t^\nu,\mu_0) \ss{\mu_0((2+\rr_t^\nu) \Phi_t'(g^\nu+c_*)f_t^\nu)} \\
&\le  k_5+k_5 t \W_2(\mu_t^\nu,\mu_0) \ss{1+\mathcal D_\nu(t)},\ \ \ t\ge t_0,\ \nu\in \scr P_0.
\end{align*}
By Young's inequality, we find $k_6\in (1,\infty)$ such that
$$
1+\mathcal D_\nu(t)\le k_6+k_6t^2\W_2(\mu_t^\nu,\mu_0)^2,\qquad t\ge t_0,\ \nu\in\scr P_0.$$
So, to prove the lower bound estimate in Theorem \ref{T1}, it remains show
\beq\label{LBB} \inf_{t\ge t_0, \nu\in \scr P_0} t\W_2(\mu_t^\nu,\mu_0)>0.\end{equation}

(d) To prove \eqref{LBB},   by \eqref{KK} we may fix $x_0\in M^\circ$ such that $\phi_0^{-1}(x_0)>4\mu(\phi_0)$. Since  $\lim_{s\to 0}Q_s\phi_0^{-1}(x_0)=\phi_0^{-1}(x_0)$, we may fix $s_*\in(0,1/2]$ such that
\beq\label{QsL}
Q_s\phi_0^{-1}(x_0)\ge3\mu(\phi_0),\qquad 0<s\le2s_*.
\end{equation}
Since  $q_s$ is the symmetric density of the Markov semigroup $Q_s$ with respect to $\mu_0$,  we have
$$
f:=q_{s_*}(x_0,\cdot)-q_{2s_*}(x_0,\cdot) = Q_{s_*/2} \big(q_{s_*/2}(x_0,\cdot)-q_{3s_*/2}(x_0,\cdot)\big),
$$
so that     \eqref{QDG-L} and  \eqref{UC} imply $ \|\nn f\|_\infty<\infty$. Moreover, by \eqref{QS}  we obtain
$$\mu_0(fq_s^\nu)= q_{s+s_*}^\nu(x_0)-q_{s+2s_*}^\nu(x_0),\ \ \ s\ge 0.$$
Moreover, \eqref{UC} implies
$$\int_{s_*}^\infty \big(|q_s^\nu(x_0)-1|+|Q_s\phi_0(x_0)^{-1}-\mu(\phi_0)|\big)\d s<\infty.$$
Combining   these  with \eqref{TTX} and \eqref{QsL},  we derive
\beg{align*}
&t\nu_0(Q_t\phi_0^{-1})\mu_0(f\tt\rr_t^\nu)\\
& =\mu(\phi_0)\int_0^\infty\left[q_{s+s_*}^\nu(x_0)-q_{s+2s_*}^\nu(x_0)\right]\d s+\int_0^\infty\left[Q_{s+s_*}\phi_0^{-1}(x_0)-Q_{s+2s_*}\phi_0^{-1}(x_0)\right]\d s\\
& =\mu(\phi_0)\int_{s_*}^{2s_*}\bigl(q_s^\nu(x_0)-1\bigr)\d s+\int_{s_*}^{2s_*}\bigl(Q_s\phi_0^{-1}(x_0)-\mu(\phi_0)\bigr)\d s\\
& =\mu(\phi_0)\int_{s_*}^{2s_*}q_s^\nu(x_0)\d s+\int_{s_*}^{2s_*}Q_s\phi_0^{-1}(x_0)\d s-2s_*\mu(\phi_0)
 \ge s_*\mu(\phi_0)>0,\qquad \nu\in\scr P_0.
\end{align*}
By this  and \eqref{QB}, there exist $c,t_0>0$ such that%we arrive at
\beq\label{FRL}
\mu_0(f \tt\rr_t^\nu)\ge \ff c t,\qquad t\ge t_0.
\end{equation}
%\beq\label{FRL}
%\inf_{t\ge t_0, \nu\in \scr P_0} \mu_0(f\tt\rr_t^\nu)>0.
%\end{equation}
Since $\mu_t^\nu=(1+\rr_t^\nu)\mu_0$,   by Kantorovich's dual   formula, we obtain
\beg{align*}
|\mu_0(f\rr_t^\nu)|=|\mu_t^\nu(f)-\mu_0(f)|\le \|\nn f\|_\infty\W_1(\mu_t^\nu,\mu_0)\le \|\nn f\|_\infty\W_2(\mu_t^\nu,\mu_0).
\end{align*} This with $ \|\nn f\|_\infty<\infty$, \eqref{ER} and \eqref{FRL} imply \eqref{LBB}, and hence finish the proof.

\end{proof}

\section{Proof of Corollary \ref{C1}}

In this part we fix $k\in [0,d]$, let $n=d-k$ and consider $\nu\in \scr P_{0,k}$. Let
$$C_{\nu}:=\inf\big\{C>0:\ C^{-1} r^k\le \nu(B(x,r))\le C r^k,\ r\in (0,1],\ x\in {\rm sup}\nu\big\},$$
and
\beq\label{Fnu}
F_\nu:=
\begin{cases}
1+\|g^\nu\|_{H^{-1}(\mu_0)}^2,&n<6,\\
\nu(\phi_0)^{-\ff{n+2}{n-2}}\int_{M^\circ}\phi_0^{\ff{n-6}{n-2}}\d\nu,&n\ge6,
\end{cases}
\end{equation} where, by \eqref{D-INF},  $\|g^\nu\|_{H^{-1}(\mu_0)}^2=\mathcal D_\nu(\infty).$  Moreover,  let
$$
r_\nu:=\sup\left\{r\in(0,1]:\nu\bigl(\{\phi_0\ge r\}\bigr)\ge\ff12\nu(M^\circ)\right\},
$$
and, for fixed $\aa\in (6,\infty)$ when $n=6$, let
\beq\label{TNU}
t_\nu:=
\begin{cases}
\left(1+\ff{C_\nu\nu(M^\circ)}{\nu(\phi_0)^2F_\nu}\right)^{\ff d{(6-n)\wedge2}},&n<6,\\
1+\left(\ff1{C_\nu \nu(\phi_0)r_\nu^5}\right)^{\alpha/6},&n=6,\\
1+\ff1{C_\nu \nu(\phi_0)r_\nu^{n-1}},&n>6.
\end{cases}
\end{equation}
Then Corollary \ref{C1} is implied by the following result.

\beg{prp}\label{T2}
Let $k\in[0,d]$ and $n=d-k$, and fix $\aa\in (6,\infty)$ when $n=6$.   Then $F_\nu<\infty$ for $\nu\in \scr P_{0,k}$, and there exists
$K \in [1,\infty)$ such that for any  $t\ge K $ and $ \nu\in \scr P_{0,k},$
\beq\label{A1}
\W_2(\mu_t^\nu,\mu_0)^2\le K  F_\nu
\begin{cases}
t^{-2},&n<6,\\
C_\nu t^{-2}\log t,&n=6,\\
C_\nu ^{\ff4{n-2}}t^{-\ff{n+2}{n-2}},&n>6.
\end{cases}
\end{equation}
Moreover, for any $t\ge K  t_\nu$ and $\nu\in \scr P_{0,k}$,
\beq\label{A2}
\W_2(\mu_t^\nu,\mu_0)^2\ge K^{-1}F_\nu
\begin{cases}
t^{-2},&n<6,\\
C_\nu ^{-1}t^{-2}\log t,&n=6,\\
C_\nu ^{-\ff4{n-2}}t^{-\ff{n+2}{n-2}},&n>6.
\end{cases}
\end{equation}
\end{prp}

 % In the following,  we prove \eqref{A1} and \eqref{A2} respectively. For simplicity, from nor on, $c $ stands for a positive constant which may change from line to line but only  depends only on $M$, $V$ and $k$, and also on $\alpha$ when $n=6$.

\beg{proof}[Proof of \eqref{A1}] By Theorem \ref{T1} and \eqref{D-L1}, it suffices to estimate $\mu_0((g^\nu+c_*)f_t^\nu)$ for $t\ge 1$.

By \eqref{V0}, \eqref{SM} and the $Q_s$-invariance of $\mu_0$, we obtain
\beq\label{DP1} \mu_0\big((g^\nu+c_*)f_t^\nu\big)= \int_0^1 \d s \int_0^\infty \mu_0\big((q_u^\nu-1)Q_s(\Phi_t'(g^\nu+c_*))\big)\d u + c_* \mu_0\big(\Phi_t'(g^\nu+c_*)\big).\end{equation}
Since $\Phi_t'(a)=\ff{ta}{t+a}\in[0,t\wedge a]$, and $\mu_0(g^\nu)=0$ due to \eqref{G-REP}, we have
\beq\label{DP2} \mu_0\big(Q_s \Phi_t'(g^\nu+c_*)\big) = \mu_0\big(  \Phi_t'(g^\nu+c_*)\big)\le \mu_0(g^\nu+c_*)= c_*,\ \ \ s\ge 0.\end{equation}
Moreover, by \eqref{QS} and the symmetry of $Q_s$ in $L^2(\mu_0)$, we obtain
$$\mu_0\big((q_u^\nu-1)Q_s\Phi_t'(g^\nu+c_*)\big)\le \mu_0\big(q_u^\nu Q_s\Phi_t'(g^\nu+c_*)\big) =\nu_0\big(Q_{u+s} \Phi_t'(g^\nu+c_*)\big),$$
and \eqref{SQI} together with \eqref{DP2} yields
$$\big|\mu_0\big((q_u^\nu-1) Q_s\Phi_t'(g^\nu+c_*)\big)\big|\le c_* c_1 \e^{-(\ll_1-\ll_0)u},\ \ \ u\ge 1,\ s\ge 0.$$
Combining these with \eqref{DP1} and \eqref{DP2}, we find $k_1\in (1,\infty)$ such that
\beg{align*}& \mu_0\big((g^\nu+c_*)f_t^\nu\big)-k_1 \le \int_0^1\d s \int_0^1\nu_0\big(Q_{u+s} \Phi_t'(g^\nu+c_*)\big)\d u\\
&= \int_0^2 \nu_0\big(Q_{r} \Phi_t'(g^\nu+c_*)\big)\d r \int_0^{1\land r}\d s\\
&= \int_0^1 r \nu_0\big(Q_{r} \Phi_t'(g^\nu+c_*)\big)\d r+  \int_1^{2} \nu_0\big(Q_{r} \Phi_t'(g^\nu+c_*)\big)\d r.\end{align*}
Noting that \eqref{UC} and \eqref{DP2} imply
$$\nu_0\big(Q_{r} \Phi_t'(g^\nu+c_*)\big) \le \big\|Q_{r} \Phi_t'(g^\nu+c_*)\big\|_\infty\le k' \mu_0\big(\Phi_t'(g^\nu+c_*)\big)\le k'c_*,\ \ r\ge 1,$$
by combining the previous estimate with Theorem \ref{T1} and \eqref{D-L1},  we find $k_2\in (1,\infty)$ such that
\beq\label{D-U}1+\mathcal D_\nu(t) \le k_2 + k_2 \int_0^1s\,\nu_0\bigl(\Phi_t'(Q_sg^\nu+c_*)\bigr)\d s,\ \ \ t\ge t_0.
\end{equation}
We will apply \eqref{QGA} to estimate the upper bound.

By the definition of $C_v$, we find $k_3\in (1,\infty)$ such that
\beg{align*}
&\int_M\e^{-\rr(x,y)^2/(c_1u)}\nu(\d y)= \sum_{j=0}^\infty \int_{\{j\ss u\le \rr(x,\cdot)<(j+1)\ss u\}} \e^{-\rr(x,\cdot)^2/(c_1u)}\nu(\d y)\\
& \le C_\nu u^{k/2}\sum_{j=0}^\infty(j+1)^k\e^{-j^2/c_1}
\le k_3 C_\nu u^{k/2},\qquad x\in S_\nu,\quad 0<u\le1.\end{align*}
Recalling that $n=d-k$,  this and the upper bound in \eqref{QGA} imply
\beq\label{qG}
\beg{split}
&q_u^\nu(x)
 =\ff1{\nu(\phi_0)}\int_{M^\circ}q_u(x,y)\phi_0(y)\nu(\d y)
 \le\ff{c_1u^{-d/2}}{\nu(\phi_0)\phi_0(x)}\int_M\e^{-\rr(x,y)^2/(c_1u)}\nu(\d y)\\
&\le\ff{c_1k_3C_\nu } {\nu(\phi_0)\phi_0(x)}u^{-n/2},\qquad x\in S_\nu,\quad 0<u\le1.
\end{split}
\end{equation}
Since $C_\nu\ge1$ and $\nu$ is a probability measure, we have
\beq\label{QQG}
\ff{C_\nu}{\nu(\phi_0)\phi_0(x)}\ge\|\phi_0\|_\infty^{-2},\qquad x\in M^\circ.\end{equation}
Thus, by \eqref{qG}, \eqref{G-REP} and \eqref{UC}, we find $k_4,k_5\in (1,\infty)$ such that
for $x\in S_\nu$ and $0<s\le1$,
\beq\label{QG}\beg{split}
&0\le Q_sg^\nu(x)+c_* =c_*+\int_s^\infty(q_u^\nu(x)-1)\d u\\
&\le k_4+\ff{k_4C_\nu}{\nu(\phi_0)\phi_0(x)}\int_s^1u^{-n/2}\d u \\
 &\le\ff{k_5C_\nu}{\nu(\phi_0)\phi_0(x)}\left(1+\int_s^1u^{-n/2}\d u\right).
\end{split}
\end{equation}
In the following,  we consider $n<6$ and $n\ge 6$ respectively.

(a) Let $n<6$.
By  \eqref{QG}  and $\Phi_t'(a)\le a$, we obtain
\beg{align*}
%&
\int_0^1s\,\nu_0(Q_sg^\nu+c_*)\d s \le\ff {k_5C_\nu\nu(M^\circ)}{\nu(\phi_0)^2}\int_0^1s\left(1+\int_s^1u^{-n/2}\d u\right)\d s<\infty.
  %\\
%&= \ff{k_5C_\nu  \nu(M^\circ)}{\nu(\phi_0)^2}\left(\ff12+\ff1{6-n}\right)<\infty.
\end{align*}
So,  \eqref{D-U} implies
$$
  F_\nu:=1+\mathcal D_\nu(\infty)=1+\lim_{t\to\infty} \mathcal D_\nu(t)<\infty,
$$
and \eqref{A1}   follows from Theorem \ref{T1} and  \eqref{D-U}.

(b) Let $n\ge6$. By $\Phi'_t(a)=\ff{ta}{t+a}$, for any $a,b>0$ we have
\beq\label{PHI-I}
\beg{split}
&\int_0^a s\,\Phi_t'\bigl(bs^{-(n-2)/2}\bigr)\d s\\
&\quad=
\begin{cases}
\ff b2\log\left(1+\ff{a^2t}b\right),&n=6,\\
b^{4/(n-2)}t^{(n-6)/(n-2)}\int_0^{a(t/b)^{2/(n-2)}}\ff{u}{1+u^{(n-2)/2}}\d u,&n>6.
\end{cases}
\end{split}
\end{equation}
We take
$$
a=1,\qquad b=\ff{k_5 C_\nu}{\nu(\phi_0)\phi_0(x)}.
$$
When  $n=6$, \eqref{QQG}  implies $\log(1+t/b)\le k_6\log t$ for some $k_6\in (1,\infty)$ and all $t\ge 2$.
For $n>6$, we have
$$
\int_0^\infty\ff{u}{1+u^{(n-2)/2}}\d u<\infty.
$$
Thus, by \eqref{QG}, \eqref{PHI-I} and the monotonicity of $\Phi_t'$, we find $k_7\in (1,\infty)$ such that for any $t\ge 2$ and $\nu\in \scr P_{0,k}$,
\beq\label{SP}
\int_0^1s\,\Phi_t'(Q_sg^\nu(x)+c_*)\d s\le k_7\cdot
\begin{cases}
\ff{C_\nu}{\nu(\phi_0)\phi_0(x)}\log t,&n=6,\\
 \left(\ff{C_\nu}{\nu(\phi_0)\phi_0(x)}\right)^{4/(n-2)}t^{(n-6)/(n-2)},&n>6.
\end{cases}
\end{equation}
By the definition of $\nu_0$, $(n-6)/(n-2)\in[0,1)$, $\phi_0$ is bounded and $\nu(\phi_0)>0$,
\beq\label{F-I}
\int_{M^\circ}(\nu(\phi_0)\phi_0(x))^{-4/(n-2)}\nu_0(\d x)=\nu(\phi_0)^{-(n+2)/(n-2)}\int_{M^\circ}\phi_0^{(n-6)/(n-2)}\d\nu=F_\nu<\infty.
\end{equation}
 Moreover, $\nu_0$ is a probability measure, so
\eqref{F-I} and
$\nu(\phi_0)\phi_0(x)\le\|\phi_0\|_\infty^2$
imply
$$
F_\nu\ge\|\phi_0\|_\infty^{-8/(n-2)}>0.
$$
Using \eqref{SP}, \eqref{F-I} and \eqref{D-U}, we find $k_8\in (1,\infty)$ such that for any $\nu\in \scr P_{0,k}$ and $t\ge 2$,
$$
1+\mathcal D_\nu(t)\le k_8 F_\nu\cdot
\begin{cases}
C_\nu\log t,&n=6,\\
C_\nu^{4/(n-2)}t^{(n-6)/(n-2)},&n>6.
\end{cases}
$$
Combining this with  \eqref{MAIN-D}, we  prove \eqref{A1} for $n\ge6$.\end{proof}

\beg{proof}[Proof of \eqref{A2}]     Since $\Phi_t'(a)=ta/(t+a)$, we have
$$
\nn\bigl(\Phi_t'(Q_sg^\nu+c_*)\bigr)=\ff{\nn Q_sg^\nu}{\bigl(1+t^{-1}(Q_sg^\nu+c_*)\bigr)^2}.
$$
Using $L_0Q_sg^\nu=1-q_s^\nu$, symmetry, the integration by parts formula, $\mu_0(Q_sg^\nu+c_*)=c_*$ and the definition of $\mathcal D_\nu(t)$ in \eqref{D-nu},
\beq\begin{split}\label{D-Q}
&\int_0^1s\,\nu_0\left(Q_s\bigl(\Phi_t'(Q_sg^\nu+c_*)\bigr)\right)\d s\\
&\quad=\int_0^1s\,\mu_0\bigl(q_s^\nu\Phi_t'(Q_sg^\nu+c_*)\bigr)\d s\\
&\quad=-\int_0^1s\,\mu_0\bigl(\Phi_t'(Q_sg^\nu+c_*)L_0Q_sg^\nu\bigr)\d s+\int_0^1s\,\mu_0\bigl(\Phi_t'(Q_sg^\nu+c_*)\bigr)\d s\\
&\quad=\int_0^1s\,\mu_0\left(\ff{|\nn Q_sg^\nu|^2}{\bigl(1+t^{-1}(Q_sg^\nu+c_*)\bigr)^2}
\right)\d s+\int_0^1s\,\mu_0\bigl(\Phi_t'(Q_sg^\nu+c_*)\bigr)\d s\\
&\quad\le\int_0^\infty s\,\mu_0\left(\ff{|\nn Q_sg^\nu|^2}
{\bigl(1+t^{-1}(Q_sg^\nu+c_*)\bigr)^2}\right)\d s+\ff12c_*\\
&\quad\le c\bigl(1+\mathcal D_\nu(t)\bigr),\qquad t\ge1.
\end{split}\end{equation}
Below we consider $n<6$ and $d\ge 6$ respectively.

(a) Let $n<6$. By \eqref{V0}, \eqref{QS} and $\pp_s Q_s=L_0Q_s$, we obtain
$$L_0Q_s g^\nu= \int_0^\infty \pp_s (q_{u+s}^\nu-1)\d s= 1-q_s^\nu,\ \ \ s>0.$$
Combining this with  the integration by parts formula, we derive
\beg{align*}
\mu_0\bigl(|\nn Q_sg^\nu|^2\bigr)&=-\mu_0\bigl(Q_sg^\nu L_0Q_sg^\nu\bigr)
=\mu_0\bigl(Q_sg^\nu(q_s^\nu-1)\bigr)\\
&=\mu_0\bigl(q_s^\nu Q_sg^\nu\bigr)=\nu_0\bigl(Q_s(Q_sg^\nu)\bigr)=\nu_0(Q_{2s}g^\nu).
\end{align*}
By \eqref{QG} for $2s$ in place of $s$, and using the definition of $\nu_0$ in \eqref{V0}, we obtain
\beq\label{D-SM}
\beg{split}
&\int_0^\vv s\,\mu_0\bigl(|\nn Q_sg^\nu|^2\bigr)\d s =\int_0^\vv s\, \nu_0(Q_{2s}g^\nu)\d s\le\int_0^\vv s\,\nu_0(Q_{2s}g^\nu+c_*)\d s\\
& \le\ff{k_5C_\nu\nu(M^\circ)}{\nu(\phi_0)^2}\int_0^\vv s\left(1+\int_s^1u^{-n/2}\d u\right)\d s
 \le\ff{k_5C_\nu \nu(M^\circ)}{\nu(\phi_0)^2}\vv^{\ff{(6-n)\wedge2}2},\qquad 0<\varepsilon\le\ff14.
\end{split}
\end{equation}
Taking
$$
\vv:=\vv_1 \left(1+\ff{A\nu(M^\circ)}{\nu(\phi_0)^2F_\nu}\right)^{-\ff2{(6-n)\wedge2}}
$$ for small enough $\vv\in (0,\ff 1 4]$,
we derive
\beq\label{QF}
4\int_0^\varepsilon s\,\mu_0\bigl(|\nn Q_sg^\nu|^2\bigr)\d s\le\ff12F_\nu.
\end{equation}
Moreover, by the $L^\infty$-contraction property of $Q_s$ and \eqref{UC}, we find $d_1,d_2\in (1,\infty)$ such that
$$
\beg{aligned}
\sup_{s\ge\varepsilon}\|Q_sg^\nu\|_\infty\le\|Q_\varepsilon g^\nu\|_\infty\le d_1+\int_\varepsilon^1\|q_u^\nu\|_\infty\d u\le d_2
\left(1+\ff{C_\nu\nu(M^\circ)}{\nu(\phi_0)^2F_\nu}\right)^{\ff d{(6-n)\wedge2}}.
\end{aligned}
$$
By the definition of $t_\nu$ in \eqref{TNU}, we find $d_3\in (1,\infty)$ such that
$$
\|Q_sg^\nu\|_\infty\le t,\qquad s\ge\varepsilon, \ \nu\in \scr P_{0,k}, \ t\ge d_3t_\nu.
$$
Combining  this with  \eqref{D-nu} and \eqref{QF},  we arrive at
\beg{align*}
&1+\mathcal D_\nu(t)
 \ge1+\int_\varepsilon^\infty s\,\mu_0\bigl(|\nn Q_sg^\nu|^2\bigr)\d s\\
&\ge1+\ff14\left[\mathcal D_\nu(\infty)-\ff12F_\nu\right]=\ff34+\ff18F_\nu\ge\ff18F_\nu,\qquad \nu\in \scr P_{0,k},\ t\ge d_3 t_\nu.
\end{align*}
This together with  \eqref{MAIN-D} implies       \eqref{A2}  for  $n<6$.

(b) Let $n\ge 6$.  Let $c_1>1$ be the constant
in \eqref{QBL}. Fix $x\in S_\nu$ and let
$$
s_x:=\ff1{8c_1^2}\bigl(1\wedge\rr_\pp(x)^2\bigr).
$$
Then for any $s\in (0,s_x]$, $y\in B(x,c_1\sqrt s)$ and $z\in B(x,\sqrt s)$,
$$\rr_\pp(y)\ge\ff12\rr_\pp(x)\ge\sqrt{2s},\qquad\rr_\pp(z)\ge\ff12\rr_\pp(x)\ge\sqrt{2s},\qquad
\rr(y,z)\le(c_1+1)\sqrt s,
$$
and  by  \eqref{KK}, we find $c_2\in (0,\infty)$ such that
$$
\phi_0(y)+\sqrt u\le c_2\phi_0(x),\qquad \phi_0(z)+\sqrt u\le c_2\phi_0(z),\ \ \ u\le 2 s.
$$
Combining these with \eqref{QGA} and the definition of $C_\nu$, we find $c_2,c_3,c_4\in (0,\infty)$ such that for $s\in (0,s_x], u\in [s,2s], y\in B\big(x,c_1\ss s\big)$,
\beg{align*}
q_u^\nu(y)
&\ge\ff1{\nu(\phi_0)}\int_{B(x,\sqrt s)}q_u(y,z)\phi_0(z)\nu(\d z)\\
&\ge\ff{u^{-d/2}}{c_1\nu(\phi_0)}\int_{B(x,\sqrt s)}\ff{\phi_0(z)\e^{- c_1\rr(y,z)^2/u}}{(\phi_0(y)+\sqrt u)(\phi_0(z)+\sqrt u)}\nu(\d z)\\
&\ge\ff{c_3s^{-d/2}}{\nu(\phi_0)\phi_0(x)}\nu(B(x,\sqrt s))\ge\ff{c_4s^{-n/2}}{C_\nu \nu(\phi_0)\phi_0(x)}.
\end{align*}
By this with \eqref{QSG} and $2s_x\le1$, there exists a constant $c_4>0$ such that
$$
\beg{aligned}
Q_sg^\nu(y)+c_*&\ge\int_s^1q_u^\nu(y)\d u\ge\int_s^{2s}q_u^\nu(y)\d u\\
&\ge\ff{c_4s^{-(n-2)/2}}{C_\nu \nu(\phi_0)\phi_0(x)},\qquad y\in B(x,c_1\sqrt s),\ s\in (0,s_x].
\end{aligned}
$$
Consequently, by \eqref{QBL} and the monotonicity of $\Phi_t'$,
\beq\label{QPL}
Q_s\bigl(\Phi_t'(Q_sg^\nu+c_*)\bigr)(x)\ge\ff12\Phi_t'\left(\ff{c_4s^{-(n-2)/2}}{C_\nu \nu(\phi_0)\phi_0(x)}\right),\ \ \
x\in S_\nu,\ s\in (0,s_x].
\end{equation}
By the definition of $r_\nu$, we have
$$
0<r_\nu\le1,\qquad\nu\bigl(\{\phi_0\ge r_\nu\}\bigr)\ge\ff12\nu(M^\circ).
$$
Since $(n-6)/(n-2)\ge0$,
\beg{align*}
&\int_{\{0<\phi_0<r_\nu\}}\phi_0^{(n-6)/(n-2)}\d\nu \le r_\nu^{(n-6)/(n-2)}\nu\bigl(\{0<\phi_0<r_\nu\}\bigr)\\
&\le r_\nu^{(n-6)/(n-2)}\nu\bigl(\{\phi_0\ge r_\nu\}\bigr)
 \le\int_{\{\phi_0\ge r_\nu\}}\phi_0^{(n-6)/(n-2)}\d\nu.
\end{align*}
Hence
\beq\label{F-L}
\int_{\{\phi_0\ge r_\nu\}}\phi_0^{(n-6)/(n-2)}\d\nu\ge\ff12\int_{M^\circ}\phi_0^{(n-6)/(n-2)}\d\nu.
\end{equation}
By \eqref{KK}, we may choose a fixed $\vv_0\in(0,1)$ such that
$$
\vv_0 r_\nu^2\le\ s_x,\qquad \phi_0(x)\ge r_\nu.
$$
Thus, by \eqref{D-Q} and \eqref{QPL}, we find $c_5\in (0,\infty)$ such that
\beq\label{D-L}
1+\mathcal D_\nu(t)\ge c_5\int_{\{\phi_0\ge r_\nu\}}\nu_0(\d x)\int_0^{\vv_0r_\nu^2}s\,\Phi_t'\left(\ff{c_4s^{-(n-2)/2}}{C_\nu \nu(\phi_0)\phi_0(x)}\right)\d s,\qquad t\ge1.
\end{equation}
Taking
$$
a=\vv_0r_\nu^2,\qquad b=\ff{c_4}{C_\nu \nu(\phi_0)\phi_0(x)}
$$
in \eqref{PHI-I}, and substituting into \eqref{D-L}, we find $c_6\in (0,\infty)$ such that
\beq\label{D-6}
\beg{split}
1+\mathcal D_\nu(t)&\ge\ff{c_6}{C_\nu \nu(\phi_0)^2}\int_{\{\phi_0\ge r_\nu\}}\log\left(1+c_6C_\nu t\nu(\phi_0)\phi_0(x)r_\nu^4\right)\d\nu,\qquad n=6,
\end{split}
\end{equation}
and for $n>6$
\beq\label{D-N}
\beg{split}
1+\mathcal D_\nu(t)&\ge\ff{c_6t^{(n-6)/(n-2)}}{C_\nu ^{4/(n-2)}\nu(\phi_0)^{(n+2)/(n-2)}}\\
&\times\int_{\{\phi_0\ge r_\nu\}}\phi_0(x)^{(n-6)/(n-2)}\left[\int_0^{c_6r_\nu^2(C_\nu t\nu(\phi_0)\phi_0(x))^{2/(n-2)}}
\ff{u}{1+u^{(n-2)/2}}\d u\right]\nu(\d x).
\end{split}
\end{equation}

When $n=6$, for any $\aa>6$, by the definition of $t_\nu$ we may choose large enough constant $K\in (1,\infty)$ such that   for $t\ge K t_\nu$, we have
$$
\log\left(1+c_6 C_\nu t\nu(\phi_0)r_\nu^5\right)\ge\log\bigl(1+t^{1-6/\alpha}\bigr)
\ge\left(1-\ff6\alpha\right)\log t.
$$
Combining this with  \eqref{D-6} and \eqref{F-L},  we find $c_7\in (0,\infty)$ such that
$$
1+\mathcal D_\nu(t)\ge c_7 C_\nu ^{-1}F_\nu\log t,\qquad t\ge K t_\nu.
$$ So, \eqref{A2} holds.

Finally, let $n>6$.   By the definition of $t_\nu$ we may choose large enough constant $K$ such that for $t\ge  Kt_\nu$
and $\phi_0(x)\ge r_\nu$ we have
$$
\int_0^{c_6r_\nu^2(C_\nu t\nu(\phi_0)\phi_0(x)/c_4)^{2/(n-2)}}\ff{u}{1+u^{(n-2)/2}}\d u\ge \int_0^1\ff{u}{1+u^{(n-2)/2}}\d u>0.
$$
Combining this with  \eqref{D-N} and \eqref{F-L}, we find a large constant $K\in (1,\infty)$ such that
$$
1+\mathcal D_\nu(t)\ge K^{-1} C_\nu ^{-4/(n-2)}F_\nu t^{(n-6)/(n-2)},\qquad  t\ge Kt_\nu,
$$
which together with   \eqref{MAIN-D}  proves \eqref{A2} for $n> 6$.
\end{proof}

\section{Proof of Corollary \ref{C2}}

\beg{proof}[Proof of Corollary \ref{C2}]
For $\nu=\dd_x$, $x\in M^\circ$, we have $\dd_x(B(x,r))=1$ for any $r>0$. Thus $\dd_x\in \scr P_{0,0}$ and
\beq\label{OP01} C_{\dd_x}=1,\ \   \ \ x\in M^\circ.\end{equation}  Moreover, we have
$\dd_x(\{\phi_0\ge r\})= 1_{\{\phi_0(x)\ge r\}},$ so that \eqref{KK} implies
\beq\label{OP02} r_{\dd_x}= \phi_0(x)\asymp \rr_\pp(x),\ \ \ x\in M^\circ. \end{equation}  By Proposition \ref{T2}, it suffices to estimate $F_{\dd_x}$.

By \eqref{KK} and the definition of $F_\nu$ in \eqref{Fnu}, it's easy to see that  when $d\ge 6$
\beq\label{OP03} F_{\dd_x}=\phi_0(x)^{-\ff 8 {d-2}}\asymp \rr_\pp(x)^{-\ff 8 {d-2}},\ \ x\in M^\circ. \end{equation} It remains to estimate $\|g^{\dd_x}\|_{H^{-1}(\mu_0)}^2$ for $d<6$. By \eqref{D-INF}, \eqref{GREEN} \eqref{V0} and \eqref{QE}, we have
\begin{equation*}\begin{split}
\|g^{\dd_x}\|_{H^{-1}(\mu_0)}^2&=4\int_0^\infty s \mu_0(|\nn Q_s G_0(x,\cdot)|^2)\d s\\
&=4\int_0^\infty s \sum_{m=1}^\infty \ff 1 {\ll_m-\ll_0} \e^{-2(\ll_m-\ll_0)s}\psi_m^2(x)\d s\\
&=\sum_{m=1}^\infty \ff 1 {(\ll_m-\ll_0)^3}\psi_m^2(x),\ \ x\in M^\circ.
\end{split}\end{equation*}
Since $$q_s(x,x)-1=\sum_{m=1}^\infty \e^{-(\ll_m-\ll_0)s}\psi_m^2(x),\  \ x\in M^\circ,\ s>0,$$
and
$$\ff 1 {(\ll_m-\ll_0)^3}=\ff 1 2 \int_0^\infty s^2\e^{-(\ll_m-\ll_0)s}\d s,\ \ m\ge 1,$$
it follows that
\beq\label{H-1}
\|g^{\dd_x}\|_{H^{-1}(\mu_0)}^2=\ff 1 2 \int_0^\infty s^2(q_s(x,x)-1)\d s,\ \ x\in M^\circ.
\end{equation}
By this with \eqref{UC} and $1\wedge \rr_\pp(x)\asymp\rr_\pp(x)$ since $\rr_\pp$ is bounded, there exist some constants $c_1,c_2,c_3>0$ such that for any $x\in M^\circ$,
\beq\begin{split}\label{gd2}
\|g^{\dd_x}\|_{H^{-1}(\mu_0)}^2&\le c_1\int_0^1 s^2 \min\{s^{-d/2}\rr_{\pp}^{-2}(x),s^{-(d+2)/2}\}\d s+c_1\\
&\le c_2\int_0^{1\wedge \rr_{\pp}^2(x)} s^{2-d/2}\rr_\pp^{-2}(x)\d s+c_2\int_{1\wedge \rr_\pp^2(x)}^1 s^{1- d/2}\d s+c_2\\
&\le c_3
\begin{cases}
1,\ \ &d\le3,\\
\log(1+\rr_\pp(x)^{-1}),\ \ &d=4,\\
\rr_\pp^{-1}(x),\ \ &d=5.
\end{cases}
\end{split}\end{equation}
By \eqref{DHL}, there exists a constant $s_0\in(0,1]$ such that
$$q_s(x,x)\ge 2,\ \ s\in(0,s_0],$$
and therefore,
$$q_s(x,x)-1\ge \ff 1 2 q_s(x,x)\ge \ff {cs^{-d/2}}{s+\rr_\pp^2(x)},\ \ s\in(0,s_0],$$
for some constant $c>0$. Substituting this into \eqref{H-1} and since $\rr$ is bounded, we can find $c_4,c_5>0$ such that for any $x\in M^\circ$,
$$\|g^{\dd_x}\|_{H^{-1}(\mu_0)}^2\ge c_4 \int_0^{s_0}
\ff{ s^{2-d/2}}{s+\rr_\pp^2(x)}\d s\ge  c_5
\begin{cases}
1,\ \ &d\le3,\\
\log(1+\rr_\pp(x)^{-1}),\ \ &d=4,\\
\rr_\pp^{-1}(x),\ \ &d=5.
\end{cases}
$$
Combining this with  \eqref{OP03} for $d\ge 6$ and \eqref{gd2},  when $\nu=\dd_x$  we deduce from \eqref{Fnu} that
$$F_{\dd_x} \asymp \begin{cases}
1,\ \ &d\le3,\\
\log(1+\rr_\pp(x)^{-1}),\ \ &d=4,\\
\rr_\pp^{-1}(x),\ \ &d=5,\\
\rr_\pp^{-\ff 8{d-2}}.\ \ &d\ge 6.
\end{cases}$$
According to Proposition \ref{T2}, this together with \eqref{OP01} and \eqref{OP02}  implies the desired estimates for some constant $c\in (1,\infty)$, which   depends on $\aa\in (6,\infty)$ only when $d=6$.
\end{proof}

\beg{thebibliography}{999}

\bibitem{AT} L. Ambrosio, F. Stra, D. Trevisan, \emph{ A PDE approach to a 2-dimensional matching problem,} Probab. Theory Related Fields 173(2019), 433-477.
\bibitem{AB15}
D. Applebaum, R. Ba\~nuelos,
\emph{Probabilistic approach to fractional integrals and the
Hardy-Littlewood-Sobolev inequality},
in: Analytic Methods in Interdisciplinary Applications,
Springer Proc. Math. Stat. 116(2015), 17-40.

\bibitem{ATW09} M. Arnaudon, A. Thalmaier, F.-Y. Wang,  \emph{Gradient estimates and Harnack inequalities on non-compact Riemannian manifolds, } Stoch. Proc. Appl. 119(2009), 3653-3670.
\bibitem{BL}  S. Bobkov, M. Ledoux, \emph{One-dimensional empirical measures, order statistics, and Kantorovich transport distances,}
Mem. Amer. Math. Soc. 261(2019), Paper No. 1259, 126pp.
\bibitem{Chavel} I. Chavel, \emph{Eigenvalues in Riemannian Geometry,} Academic Press, 1984.

\bibitem{FS} A. Ferreira, I. Salavessa, \emph{Dirichlet principal eigenvalue comparison theorems in geometry with torsion,}
J. Math. Anal. Appl. 453(2017),  700-723.
\bibitem{FG} N. Fournier, A. Guillin, \emph{On the rate of convergence in Wasserstein distance of the empirical measure,}
Probab. Theory Relat. Fields 162(2015), 707-738.
%\bibitem{Bayle} V. Bayle, \emph{Propri\'{e}t\'{e}s de concavit\'{e} du profil isop\'{e}rim\'{e}trique et applications,} PhD thesis, Institut Fourier (Grenoble), 2003.

%\bibitem{CHT}
%Y. Colin de Verdi\`ere, L. Hillairet and E. Tr\'elat,
%\emph{Small-time asymptotics of hypoelliptic heat kernels near the
%diagonal, nilpotentization and related results},
%Ann. Henri Lebesgue 4(2021), 897--971.

%\bibitem{G06}
%A. Grigor'yan,
%\emph{Heat kernels on weighted manifolds and applications},
%in: The Ubiquitous Heat Kernel, Contemp. Math. 398,
%Amer. Math. Soc., Providence, RI, 2006, 93--191.
%\bibitem{Kendall}  W. S. Kendall,
%\emph{The radial part of Brownian motion on a manifold: a
%semimartingale property,} Ann.  Probab. 15(1987), 1491-1500.
%\bibitem{Lee}
%J. M. Lee,
%\emph{Introduction to Riemannian Manifolds},
%2nd ed., Grad. Texts in Math. 176, Springer, Cham, 2018.

\bibitem{CJ}
J. Chen, S. Jian,
\emph{A remark on quasi-ergodicity of ultracontractive Markov processes}, Stat. Probab. Lett. 87(2014) 184-190.

\bibitem{HGT}  M. Huesmann, M. Goldman, D. Trevisan,    \emph{Asymptotics for random quadratic transportation costs,} arXiv:2409.08612.
\bibitem{HM} M. Huesmann, F. Mattesini,  D. Trevisan, \emph{Wasserstein asympototics for the empirical measure of fractional Brownian motion on a flat torus,} Stoch. Proc. Appl. 155(2023), 1-26.

%\bibitem{Kloeckner} B. Kloeckner,
%Approximation by finitely supported measures,
%ESAIM Control Optim. Calc. Var. 18 (2012), 343--359.

\bibitem{Ledoux}
 M. Ledoux,
\emph{On optimal matching of Gaussian samples},
J. Math. Sci. (N.Y.) 238 (2019),  495--522.
 \bibitem{LZ} M. Ledoux, J.-X. Zhu, \emph{ On optimal matching of Gaussian samples III,}
Probab. Math. Statist. 41(2021),  237-265.
\bibitem{LW} H. Li, B.  Wu, \emph{ Wasserstein convergence rates for empirical measures of subordinated processes on noncompact manifolds,} J. Theor. Probab. 36(2023), 1243-1268.
\bibitem{LW24}  H. Li, B. Wu, \emph{Wasserstein convergence for conditional empirical measures of subordinated Dirichlet diffusions on Riemannian manifolds,}
Commun. Pure Appl. Anal. 23(2024),   546-576.
\bibitem{LW26}  H. Li, B. Wu, \emph{Wasserstein convergence for empirical measures of subordinated fractional Brownian motions on the flat torus,}
 Bull. Sci. Math. 213(2026), Paper No. 103891, 44pp.

\bibitem{MS}
P. Mattila and P. Saaranen,
Ahlfors--David regular sets and bilipschitz maps,
Ann. Acad. Sci. Fenn. Math. 34 (2009), 487--502.

\bibitem{RZ}
J. Rataj and M. Z\"ahle,
\textit{Curvature Measures of Singular Sets},
Springer Monographs in Mathematics,
Springer, Cham, 2019.
\bibitem{TW98} A. Thalmaier, F.-Y. Wang,  \emph{Gradient estimates for harmonic functions on regular domains in Riemannian manifolds,} J. Funct. Anal.  155:1(1998),  109-124.
\bibitem{TWZ} D. Trevisan, F.-Y. Wang, J.-X. Zhu, \emph{Wasserstein asymptotics for empirical measures of diffusions on four dimensional closed manifolds, }  Electron. Commun. Probab. 30(2025), Paper No. 68, 13pp.

\bibitem{W21}
F.-Y. Wang,
\emph{Precise limit in Wasserstein distance for conditional empirical
measures of Dirichlet diffusion processes},
J. Funct. Anal. 280(2021), Paper No. 108998, 23pp.

\bibitem{W22}
F.-Y. Wang, \emph{Wasserstein convergence rate for empirical measures on noncompact manifolds},
Stoch. Proc. Appl. 144(2022), 271-287.

\bibitem{WJEMS}
F.-Y. Wang,
\emph{Convergence in Wasserstein distance for empirical measures of
Dirichlet diffusion processes on manifolds},
J. Eur. Math. Soc. 25(2023), 3695-3725.
\bibitem{W25} F.-Y. Wang, \emph{Wasserstein convergence rate for empirical measures of Markov processes,}     Appl. Math. Opt.  92(2025),   Paper No. 4, 41pp.
  \bibitem{W}    F.-Y. Wang, \emph{Asymptotics in Wasserstein distance for empirical measures of Markov processes,}    Comm. Math. Anal. Appl. 4(2025), 495-524.
  \bibitem{W26}  F.-Y. Wang \emph{Convergence in Wasserstein distance for empirical measures of non-symmetric subordinated diffusion processes,}  Stoch. Proc. Appl. 198(2026),
  Paper No. 104969, 36pp.
\bibitem{WW} F.-Y. Wang, B. Wu, \emph{Wasserstein convergence for empirical measures of subordinated diffusion processes on Riemannian manifolds,} Pot. Anal. 59(2023), 933-954.
\bibitem{WWZ} F.-Y. Wang, B. Wu,  J.-X. Zhu, \emph{Sharp $L^q$-convergence Rate in $p$-Wasserstein distance for empirical measures of diffusion processes,} Stoch. Proc. Appl. 195(2026), Paper No. 104869, 17pp.

\bibitem{WZ}
F.-Y. Wang, J.-X. Zhu,
\emph{Limit theorems in Wasserstein distance for empirical measures of
diffusion processes on Riemannian manifolds},
Ann. Inst. Henri Poincar\'e Probab. Stat. 59(2023), 437-475.

\bibitem{WZ25} B. Wu, J.-X. Zhu, \emph{Wasserstein convergence rates for empirical measures of random subsequence of $\{n\aa\}$,}
 Stoch. Proc. Appl. 181(2025), Paper No. 104534, 20pp.
\bibitem{Zhang} Q. S. Zhang, \emph{ The global behavior of heat kernels in exterior domains,}
J. Funct. Anal. 200(2003),   160-176.
\end{thebibliography}

\end{document}